\documentclass[11pt,a4paper,reqno]{amsart}
\usepackage{amssymb,amsmath,amsthm}
\usepackage[utf8]{inputenc}
\usepackage[T1]{fontenc}
\usepackage{enumerate}
\usepackage[all]{xy}
\usepackage{fullpage}
\usepackage{comment}
\usepackage{array}
\usepackage{longtable}
\usepackage{stmaryrd}
\usepackage{mathrsfs} 
\usepackage{xcolor}
\usepackage{mathtools}

\def\dim{\mathrm{dim}}

\usepackage{hyperref}
\hypersetup{
    colorlinks=true,
    linkcolor=blue,
    filecolor=red,  
    citecolor=green,
    urlcolor=cyan,
    pdftitle={Orthogonal Sets},
    pdfpagemode=FullScreen,
    }

\usepackage{cleveref}
\crefformat{section}{\S#2#1#3}
\crefformat{subsection}{\S#2#1#3}
\crefformat{subsubsection}{\S#2#1#3}
\usepackage{enumitem}
\usepackage{tikz}

\newtheorem*{claim*}{Claim}

\newtheorem{theorem}{Theorem}[section]
\newtheorem{lemma}[theorem]{Lemma}
\newtheorem{corollary}[theorem]{Corollary}
\newtheorem{proposition}[theorem]{Proposition}
\newtheorem{example}[theorem]{Example}

\theoremstyle{definition}
\newtheorem{definition}[theorem]{Definition}

\theoremstyle{remark}
\newtheorem{remark}[theorem]{Remark}
\newtheorem{question}[theorem]{Question}

\makeatletter
\newcommand{\subalign}[1]{%
  \vcenter{%
    \Let@ \restore@math@cr \default@tag
    \baselineskip\fontdimen10 \scriptfont\tw@
    \advance\baselineskip\fontdimen12 \scriptfont\tw@
    \lineskip\thr@@\fontdimen8 \scriptfont\thr@@
    \lineskiplimit\lineskip
    \ialign{\hfil$\m@th\scriptstyle##$&$\m@th\scriptstyle{}##$\hfil\crcr
      #1\crcr
    }%
  }%
}
\makeatother

\numberwithin{equation}{section}

\title{On the maximum size of ultrametric orthogonal sets over discrete valued fields}

\author{Noy Soffer Aranov}
\email{noyso@campus.technion.ac.il}
\address{Department of Mathematics, Technion - Israel Institute of Technology, Haifa, Israel}

\author{Angelot Behajaina}
\email{abeha@campus.technion.ac.il}
\address{Department of Mathematics, Technion - Israel Institute of Technology, Haifa, Israel}
\address{Department of Mathematics and Computer Science, The Open University of Israel}
\begin{document}

\maketitle
\begin{abstract}
Let $\mathcal{K}$ be a discrete valued field with finite residue field. In analogy with orthogonality in the Euclidean space $\mathbb{R}^n$, there is a well-studied notion of "ultrametric orthogonality" in $\mathcal{K}^n$. In this paper, motivated by a question of Erd{\H{o}}s in the real case, given integers $k \geq  \ell \geq 2$, we investigate the maximum size of a subset $S \subseteq \mathcal{K}^n \setminus\{{\bf 0}\}$ satisfying the following property: for any $E \subseteq S$ of size $k$, there exists $F \subseteq E$ of size $\ell$ such that any two distinct vectors in $F$ are orthogonal. Other variants of this property are also studied.
\end{abstract}
\section{Introduction}
Let $\mathbb{R}^n$ ($n \geq 2$) be the Euclidean space equipped with the classical scalar product $\mathbf{u}\cdot \mathbf{v}=\sum_{i=1}^nu_iv_i$; the vectors $\mathbf{u}$ and $\mathbf{v}$ are \emph{orthogonal} if $\mathbf{u}\cdot \mathbf{v}=0$. It is well-known that the maximum size of a set $S \subseteq \mathbb{R}^n \setminus\{{\bf 0}\}$, such that any two distinct vectors of $S$ are orthogonal is $n$. More generally, Erd{\H{o}}s asked the following question.  
\begin{question}[Erd{\H{o}}s]\label{quest:erd}
Let $k \geq \ell \geq 2$. We let $\alpha_n^{(k,\ell)}$ be the maximum size of a subset $S \subseteq \mathbb{R}^n \setminus \{{\bf 0}\}$ satisfying the following property: for any $E \subseteq S$ of size $k$, there exists $F \subseteq E$ of size $\ell$ such that any two distinct vectors in $F$ are orthogonal. What is $\alpha_n^{(k,\ell)}$?
\end{question}
\noindent
This question has attracted a lot of attention in the past decades (see for example \cite{AS99,Dea11,FS92,Ros95}). Note that the observation at the beginning can be reformulated as $\alpha_n^{(2,2)}=n$. In \cite{FS92}, F\"{u}redi and Stanley proved that $\alpha_2^{(k,2)}=2k-2$. Furthermore, Rosenfeld showed \cite{Ros95} that $\alpha_n^{(3,2)}=2n$, a result which was reproved by Deaett \cite{Dea11}.

Given a finite field $\mathbb{F}_q$ where $q$ is an odd prime power, Ahmadi and Mohammadian \cite{AM16} investigated the analogue of the problem of Erd{\H{o}}s in $\mathbb{F}_q^n$ ($n \geq 2$), equipped with a symmetric non-degenerate bilinear form. In this setting, they computed $\alpha_n^{(2,2)}$. Moreover, they provided an upper bound on $\alpha_n^{(3,2)}$. In \cite{MP22}, Mohammadian and Petridis improved this upper bound.

A positive characteristic analogue of $\mathbb{R}$ that is well-studied in geometry of numbers is the field $\mathbb{F}_q((x^{-1}))$ of Laurent series in $x^{-1}$ over a finite field $\mathbb{F}_q$ (see \cite{KST}, \cite{Mah}, \cite{RW}, \cite{Arm}, \cite{Ara}). However, orthogonality with respect to bilinear forms on $\mathbb{F}_q((x^{-1}))^n$ are different. In this work, we consider the  "ultrametric orthogonality" (see Definition \ref{def:orthvect}), which makes sense in any ultrametric discrete valued field. This is inspired by the following property: in $\mathbb{R}^n$ equipped with the $\mathcal{L}_2$-norm, $\mathbf{u},\mathbf{v}\in \mathbb{R}^n$ are orthogonal if and only if, for every $a,b\in \mathbb{R}$, we have $\Vert a\mathbf{u}+b\mathbf{v}\Vert^2_2=\Vert a\mathbf{u}\Vert^2_2+\Vert b\mathbf{v}\Vert^2_2$.

Now let $\mathcal{K}$ be a discrete valued field with valuation $\nu:\mathcal{K} \rightarrow \mathbb{Z} \cup \{0\}$. The sets 
$$\mathfrak{m}=\{ \lambda \in \mathcal{K} \mid \nu(\lambda)>0\},
$$
$$\mathcal{O}=\{ \lambda \in \mathcal{K} \mid \nu(\lambda) \geq 0\},
$$ and  
$$\kappa=\mathcal{O}/\mathfrak{m}$$ denote, respectively, the maximal ideal, the valuation ring and the residue field of $\mathcal{K}$. In the following, we assume that $\kappa$ is {\bf finite}; for that let us suppose $\kappa=\mathbb{F}_q$ where $q$ is a prime power. Here are some examples of discrete valued fields. 
\begin{example}\label{ex:dvr}
\begin{itemize}
    \item Let $\mathcal{K}=\mathbb{F}_q((x^{-1}))$ ($q$ a prime power) be the field of Laurent series in $x^{-1}$ and let $\nu$ be the valuation defined by $\nu(\sum_{i}a_ix^{-i})=\min(i \mid a_i\neq 0)$. In this case, $\mathcal{O}=\mathbb{F}_q[[x^{-1}]]$, $\mathfrak{m}=x^{-1} \mathcal{O}$ and $\kappa=\mathbb{F}_q$. 
    \item Let $\mathcal{K}=\mathbb{Q}_p$ ($p$ a prime number) be the field of $p$-adic numbers and let $\nu$ be the $p$-adic valuation. In this case, $\mathcal{O}=\mathbb{Z}_p$, $\mathfrak{m}=p \mathbb{Z}_p$ and $\kappa=\mathbb{F}_p$.
\end{itemize}
\end{example}
On $\mathcal{K}$, we consider the utrametric absolute value defined by:
$$
\left\vert \lambda\right\vert =q^{-\nu(\lambda)},\,\,\textrm{for all}\,\,\lambda \in \mathcal{K}.
$$ This gives rise to the \emph{infinity norm} on $\mathcal{K}^n$ ($n \geq 1$) given by:
$$
\Vert{\bf v}\Vert=\max_{i \in [n]}\vert v_i\vert ,\,\, \textrm{for all}\,\, {\bf v}=(v_1,\dots,v_n) \in \mathcal{K}^n.
$$  Recall:
\begin{lemma}[Ultrametric Inequality]
\label{lem:UMIneq}
\begin{enumerate}
    \item For any $\alpha_1,\dots,\alpha_{\ell}\in \mathcal{K}$, we have
    \begin{equation}\label{eq:ultrscal}\vert \alpha_1+\dots+\alpha_{\ell}\vert\leq \max_{i \in [\ell]}\vert \alpha_i\vert,
    \end{equation} where $[\ell]=\{1,\dots,\ell\}$. Moreover, if there exists $i_0$ such that $\vert \alpha_{i}\vert < \vert\alpha_{i_0}\vert$ for all $i \neq i_0$, then the inequality in \eqref{eq:ultrscal} is an equality. 
    \item For any $\mathbf{u}_1,\dots ,\mathbf{u}_{\ell}\in \mathcal{K}^n$, we have
    \begin{equation}\label{eq:ultrvect}
    \Vert \mathbf{u}_1+\dots+\mathbf{u}_{\ell}\Vert\leq \max_{i \in [\ell]}\Vert \mathbf{u}_i\Vert.
    \end{equation} Moreover, if there exists $i_0$ such that $\Vert {\bf u}_{i}\Vert < \Vert {\bf u}_{i_0} \Vert$ for all $i \neq i_0$, then the inequality in \eqref{eq:ultrvect} is an equality. 
\end{enumerate}
\end{lemma}
In the sequel, for two vector spaces $U$ and $V$, we denote $U\leq V$, whenever $U$ is a subspace of $V$. The "ultrametric orthogonality" we consider is the following (see \cite{KST} and \cite{RW}).
\begin{definition}\label{def:orthvect} 
\begin{itemize}
    \item The vectors $\mathbf{u}_1,\dots,\mathbf{u}_\ell \in \mathcal{K}^n \setminus \{{\bf 0}\}$ are \emph{orthogonal} if 
    $$\Vert \lambda_\ell \mathbf{u}_1+\dots + \lambda_\ell \mathbf{u}_\ell\Vert=\max_{i \in [\ell]}\Vert \lambda_i {\bf u}_i \Vert,$$
    for all $\lambda_1,\dots,\lambda_\ell \in \mathcal{K}$. In other words, every linear combination of $\mathbf{u}_1,\dots,\mathbf{u}_\ell$ satisfies the equality case of the ultrametric inequality.
    \item The subspaces $U,V \leq \mathcal{K}^n$ are \emph{orthogonal} if $\mathbf{u}$ and $\mathbf{v}$ are orthogonal, for all ${\bf u}\in U$ and ${\bf v}\in V$.
    \item The subspaces $U, V \subseteq \mathcal{K}^n$ are \emph{feebly orthogonal} if:
    \begin{itemize}
    \item either there exists ${\bf u} \in U \setminus \{{\bf 0}\}$ such that $\mathcal{K}\cdot{\bf u}$ and $V$ are orthogonal;
    \item or there exists ${\bf v} \in V \setminus \{{\bf 0}\}$ such that $\mathcal{K}\cdot{\bf v}$ and $U$ are orthogonal.
    \end{itemize}
\end{itemize}
\end{definition}
\begin{remark}
    We can view feeble orthogonality as a weaker variant of orthogonality, in which one of the subspaces contains a vector which is orthogonal to all of the vectors in the other subspace. 
\end{remark}
\noindent
The \emph{unit sphere} in $\mathcal{K}^n$ is defined by $\mathbb{B}_n=\{{\bf v}\in \mathcal{K}^n \mid \Vert{\bf v}\Vert=1\}$, and the \emph{Grassmannian} ${\rm Gr}_{s,n}(\mathcal{K})$ ($1 \leq s \leq n$) denotes the set of all $s$-dimensional subspaces of $\mathcal{K}^n$.
\begin{definition}\label{def:pairortho} 
\begin{itemize}
\item A subset $S\subseteq \mathcal{K}^n \setminus \{{\bf 0}\}$ is \emph{weakly orthogonal} if any two distinct vectors in $S$ are orthogonal. 
\item A subset $S=\{{\bf u}_1,\dots,{\bf u}_\ell\}\subseteq \mathcal{K}^n \setminus \{{\bf 0}\}$ is \emph{ orthogonal} if the vectors ${\bf u}_1,\dots,{\bf u}_\ell$ are orthogonal in the sense of Definition \ref{def:orthvect}. 
\item Let $s \geq 1$. A family $\mathcal{F} \subseteq {\rm Gr}_{s,n}(\mathcal{K})$ is  \emph{feebly orthogonal} if any two distinct subspaces in $\mathcal{F}$ are feebly orthogonal.
\end{itemize}
\end{definition}
We can generalize Definition \ref{def:pairortho}.
\begin{definition}\label{def:ortsetgdef}
Let $ k \geq \ell \geq 2$. 
\begin{itemize}
\item A subset $S \subseteq \mathcal{K}^n \setminus \{{\bf 0}\}$ is \emph{$(k,\ell)$-weakly orthogonal} if any $E \subseteq S$ of size $k$ contains a weakly orthogonal subset $F \subseteq E$ of size $\ell$. We denote by
$\Delta_{n}^{(k,\ell)}$ the maximum size of a $(k,\ell)$-weakly orthogonal set in $\mathcal{K}^n \setminus \{{\bf 0}\}$, that is,
\begin{align*}
\Delta_{n}^{(k,\ell)}&=\max\{\vert S\vert \mid S \subseteq \mathcal{K}^n\setminus\{{\bf 0}\} \,\,\textrm{is}\,\, (k,\ell)\textrm{-weakly orthogonal}\}\\
&=\max\{\vert S\vert  \mid S \subseteq \mathbb{B}_n \,\,\textrm{is}\,\, (k,\ell)\textrm{-weakly orthogonal}\}.
\end{align*}
\item A family $\mathcal{F} \subseteq {\rm Gr}_{s,n}(\mathcal{K})$ ($1 \leq s \leq n$) is \emph{$(k,\ell)$-feebly orthogonal} if any $\mathcal{E} \subseteq \mathcal{F}$ of size $k$ contains a feebly orthogonal subset $\mathcal{R} \subseteq \mathcal{E}$ of size $\ell$. We let
\begin{align*}
\Omega_{s,n}^{(k,\ell)}&=\max\{\vert \mathcal{F}\vert  \mid \mathcal{F} \subseteq {\rm Gr}_{s,n}(\mathcal{K})\,\,\textrm{is}\,\, (k,\ell)\textrm{-feebly orthogonal}\}.
\end{align*}
\end{itemize}
\end{definition}
\begin{remark}
Note that orthogonality is invariant under scalar multiplication. Hence, we restrict to vectors in $\mathbb{B}_n$ in Definition \ref{def:ortsetgdef}.
\end{remark}
\noindent
In analogy with Question \ref{quest:erd}, it is natural to ask the following.
\begin{question}\label{qst:alphadefiniqest}
Let $s \geq 1$ and let $k \geq \ell \geq 2$. What are $\Delta_{n}^{(k,\ell)}$ and $\Omega_{s,n}^{(k,\ell)}$?
\end{question}
\noindent
The two theorems below answer Question \ref{qst:alphadefiniqest}.
\begin{theorem}[Weakly orthogonal sets]\label{thm:main1} Let $k \geq \ell \geq 2$. Then
        \begin{equation}\label{eq:boundalphakln}
            \left\lfloor \frac{k-1}{\ell-1}\right\rfloor\frac{q^n-1}{q-1}\leq \Delta_{n}^{(k,\ell)}\leq \left\lfloor\frac{k-1}{\ell-1}\cdot \frac{q^n-1}{q-1}\right\rfloor;
        \end{equation} more precisely,
        \begin{equation}\label{eq:optimform}
\Delta^{(k,\ell)}_{n}=\max\left\lbrace\left.\sum_{i=1}^{\frac{q^n-1}{q-1}}t_i\,\, \right\vert \sum_{i \in I}t_i \leq k-1,\,\,\textrm{for all}\,\, I \subseteq \left[\frac{q^n-1}{q-1}\right]\,\, \textrm{of size}\,\, \ell-1\right\rbrace.
\end{equation}
As a consequence, if $\ell-1$ divides $k-1$, then $\Delta_{n}^{(k,\ell)}=\frac{k-1}{\ell-1}\cdot\frac{q^n-1}{q-1}$. In particular, 
$$\Delta_{n}^{(k,2)}=(k-1)\frac{q^n-1}{q-1}.
$$
\end{theorem}
\begin{corollary}
Let $\ell \geq 2$. When $k \rightarrow \infty$, we have $\Delta_n^{(k,\ell)}\sim \frac{1}{\ell-1}\frac{q^n-1}{q-1} \cdot k$.    
\end{corollary}
\begin{theorem}[Feebly orthogonal spaces]\label{thm:main1ss} Let $s \geq 1$ and let $k \geq \ell \geq 2$. Then
        \begin{equation}\label{eq:boundalphaklnss}
            \left\lfloor \frac{k-1}{\ell-1}\right\rfloor{n \brack s}_{q}\leq \Omega_{s,n}^{(k,\ell)}\leq \left\lfloor\frac{k-1}{\ell-1}\cdot {n \brack s}_{q}\right\rfloor,
        \end{equation} 
where ${n \brack s}_{q}=\frac{(q^{n-s+1}-1)\cdots(q^2-1)(q-1)}{(q^s-1)\cdots(q^2-1)(q-1)}$; 
more precisely,
        \begin{equation}\label{eq:optimformss}
\Omega^{(k,\ell)}_{s,n}=\max\left\lbrace\left.\sum_{i=1}^{{n \brack s}_{q}}t_i\,\,\right\vert \sum_{i \in I}t_i \leq k-1,\,\,\textrm{for all}\,\, I \subseteq \left[{n \brack s}_{q}\right]\,\, \textrm{of size}\,\, \ell-1\right\rbrace.
\end{equation}
As a consequence, if $\ell-1$ divides $k-1$, then $\Omega_{s,n}^{(k,\ell)}=\frac{k-1}{\ell-1}\cdot{n \brack s}_{q}$. In particular, $\Omega_{s,n}^{(k,2)}=(k-1){n \brack s}_{q}$.
\end{theorem}
\begin{corollary}
Let $s \geq 1$ and let $\ell \geq 2$. When $k \rightarrow \infty$, we have $\Omega_{s,n}^{(k,\ell)}\sim \frac{1}{\ell-1}{n \brack s}_q \cdot k$.    
\end{corollary}
Motivated by Definition \ref{def:ortsetgdef}, we also consider: 
\begin{definition}\label{ref:strongorthsetdef}
    Let $k \geq \ell \geq 2$. A subset $S\subseteq \mathcal{K}^n\setminus \{{\bf 0}\}$ is \emph{$(k,\ell)$-orthogonal} if any $E\subseteq S$ of size $k$ contains an orthogonal subset $F \subseteq E$ of size $\ell$. We let
    \begin{align*}
    \Theta_{n}^{(k,\ell)}&=\max\{\vert S\vert  \mid S \subseteq \mathcal{K}^n \setminus\{{\bf 0}\}\,\,\textrm{is}\,\,(k,\ell)\textrm{-orthogonal}\}\\
    &=\max\{\vert S\vert  \mid S \subseteq \mathbb{B}_n\,\,\textrm{is}\,\,(k,\ell)\textrm{-orthogonal}\}.
    \end{align*}
\end{definition}
\begin{question}\label{qst:alphabardefiniqest}
Let $ k \geq \ell \geq 2$. What is $\Theta_{n}^{(k,\ell)}$? 
\end{question}
\noindent
In order to estimate $\Theta_{n}^{(k,\ell)}$, our approach leads to analyze "independent sets" in $\mathbb{F}_q^n$; see \cite{DMM,DMMS,H,Tal,TV09}.
\begin{definition} Let $1 \leq \ell \leq k \leq q^n-1$ with $\ell \leq n$. A subset $S\subseteq \mathbb{F}_q^n \setminus \{{\bf 0}\}$ is \emph{$(k,\ell)$-independent} if every subset $X\subseteq S$ of size $k$ contains a linearly independent subset of size $\ell$, that is $\dim(\langle S \rangle) \geq \ell$. We let
$$
{\rm Ind}_q(n,k,\ell)=\max\{\vert S\vert  \mid S \subseteq \mathbb{F}_q^n\setminus\{{\bf 0}\}\,\,\textrm{is}\,\,(k,\ell)\textrm{-independent}\}.
$$
\end{definition}
\noindent
Towards Question \ref{qst:alphabardefiniqest}, we obtain partial results:
\begin{theorem}[Orthogonal sets]\label{thm:main2} Let $k \geq \ell \geq 2$.
\begin{enumerate}
\item  
    $$ \Theta_{n}^{(k,\ell)} \leq \Delta_{n}^{(k,\ell)} \leq  \left\lfloor\frac{k-1}{\ell-1}\cdot \frac{q^n-1}{q-1}\right\rfloor .$$  \label{thm:main21} 
\item  $\operatorname{Ind}_{q}(n,k,\ell)\leq\Theta_{n}^{(k,\ell)}$, with equality when $k=\ell$. Moreover, if the equality holds, then $k\leq q^{\ell-2}+1$. As a consequence ${\rm Ind}_{q}(n,k,2)=(k-1)\frac{q^n-1}{q-1}$, for all $2\leq k \leq q$. \label{thm:main22}
\item  $\left\lfloor\frac{k-1}{\ell-1}\right\rfloor\Theta_{n}^{(\ell,\ell)} \leq \Theta_{n}^{(k,\ell)}\leq (k-\ell+1)\operatorname{Ind}_{q}(n,k,\ell)$. \label{thm:main23}
\item 
$$
  \frac{q^n-1}{q^{\ell-1}-1} \leq   \limsup_{k\rightarrow \infty}\frac{\Theta_{n}^{(k,\ell)}}{k} \leq q^{n}-1.
$$ \label{thm:main24}
\item $\Theta_{n}^{(\ell,\ell)}<\Delta_{n}^{(\ell,\ell)}$. \label{thm:main25} 
\end{enumerate}
\end{theorem}
An interesting question is whether a $(k,\ell)$-orthogonal set must be weakly orthogonal. This motivates the following definition.
\begin{definition}
    Let $k \geq \ell \geq 2$. 
    A weakly orthogonal set $S\subseteq \mathcal{K}^n\setminus \{{\bf 0}\}$ is $(k,\ell)$\emph{-strongly orthogonal} if, for any set $E\subseteq S$ of size $k$, there exists an orthogonal set $F \subseteq E$ of size $\ell$. We let \begin{align*}\Gamma_{n}^{(k,\ell)}&=\max\{\vert S\vert  \mid S \subseteq \mathcal{K}^n \setminus\{{\bf 0}\}\,\,\textrm{is}\,\,(k,\ell)\textrm{-strongly orthogonal}\}\\
    &=\max\{\vert S\vert  \mid S \subseteq \mathbb{B}_n\,\,\textrm{is}\,\,(k,\ell)\textrm{-strongly orthogonal}\}.
    \end{align*}
\end{definition}
\begin{question}\label{qst:betadefiniqest}
Let $ k \geq \ell \geq 2$. What is $\Gamma_{n}^{(k,\ell)}$? 
\end{question}
\noindent
To investigate $\Gamma_{n}^{(k,\ell)}$, we also need to study "independent sets" in $\mathbb{P}^{n-1}(\mathbb{F}_q)=(\mathbb{F}_q^n\setminus\{{\bf 0}\})/\mathbb{F}_q^*$. 
\begin{definition} Let $1 \leq \ell \leq k \leq \frac{q^n-1}{q-1}$ with $\ell \leq n$. A subset $S\subseteq \mathbb{P}^{n-1}(\mathbb{F}_q)$ is \emph{$(k,\ell)$-pro-independent} if any subset $X\subseteq S$ of size $k$ contains a linearly independent subset of size $\ell$, that is $\dim(\langle S \rangle) \geq \ell$. We let
$$
{\rm Ind}_q^{\operatorname{pro}}(n,k,\ell)=\max\{\vert S\vert  \mid S \subseteq \mathbb{P}^{n-1}(\mathbb{F}_q)\,\,\textrm{is}\,\,(k,\ell)\textrm{-pro-independent}\}.
$$
\end{definition}
\begin{remark}
\label{rmk:Ind=Indpro,q=2}
    When $q=2$, we have $\mathbb{P}^{n-1}(\mathbb{F}_2)=\mathbb{F}_2^n$, so $\operatorname{Ind}_2^{\operatorname{pro}}(n,k,\ell)=\operatorname{Ind}_2(n,k,\ell)$. 
\end{remark}
\noindent
Concerning Question \ref{qst:betadefiniqest}, we get:
\begin{theorem}[Strongly orthogonal sets]\label{thm:main3}
Let $k \geq \ell \geq 2$. Then $\Gamma_{n}^{(k,\ell)}={\rm Ind}^{{\rm pro}}_q(n,k,\ell)$.
\end{theorem}
Of independent interest, we also investigate some questions regarding independent sets in $\mathbb{F}_q^n$ and $\mathbb{P}^{n-1}(\mathbb{F}_q)$:
\begin{theorem}\label{thm:mainindaffproj}
\begin{enumerate}
\item {\bf Independent sets in $\mathbb{F}_q^n$.} Let $2 \leq \ell \leq k \leq q^n-1$ with $\ell \leq n$. Then:
\begin{enumerate}
\item  $\operatorname{Ind}_q(n,k,\ell)=q^{n}-1$ if and only if $k\geq q^{\ell-1}$. In other words, $\mathbb{F}_q^{n}\setminus \{{\bf  0}\}$ is $(k,\ell)$-independent if and only if $k\geq q^{\ell-1}$;\label{thm:mainindaffproj1a}
\item the sequence $\{{\rm Ind}_q(n,k,\ell)\}_{k=\ell}^{q^{\ell-1}}$ is strictly increasing in $k$; \label{thm:mainindaffproj1b}
\item 
$$
q^n-1={\rm Ind}_q(n,k,1) =\dots={\rm Ind}_q(n,k,\lfloor \log_q(k) \rfloor +1)> \dots > {\rm Ind}_q(n,k,k)\geq n+1.
$$ \label{thm:mainindaffproj1c}
\end{enumerate}
\item {\bf Independent sets in $\mathbb{P}^{n-1}(\mathbb{F}_q)$.} Let $2 \leq \ell \leq k \leq \frac{q^n-1}{q-1}$ with $\ell \leq n$. Then:
\begin{enumerate}
\item $\operatorname{Ind}_q^{\operatorname{pro}}(n,k,\ell)=\frac{q^n-1}{q-1}$ if and only if $k\geq \frac{q^{\ell-1}-1}{q-1}$. In other words, $\mathbb{P}^{n-1}(\mathbb{F}_q)$ is $(k,\ell)$-pro-independent if and only if $k\geq \frac{q^{\ell-1}-1}{q-1}$; \label{thm:mainindaffproj2a}
\item the sequence $\{\operatorname{Ind}_q^{\operatorname{pro}}(n,k,\ell)\}_{k=\ell}^{\frac{q^{\ell-1}-1}{q-1}}$ is strictly increasing in $k$; \label{thm:mainindaffproj2b}
\item moreover,
$$
\frac{q^n-1}{q-1}={\rm Ind}^{\operatorname{pro}}_q(n,k,1)= \dots={\rm Ind}^{{\rm pro}}_{q}(n,k,\lfloor\log_q((q-1)k+1)\rfloor+1)> \dots > {\rm Ind}^{\operatorname{pro}}_q(n,k,k)\geq n+1.
$$ \label{thm:mainindaffproj2c}
\end{enumerate}
\item {\bf Independent sets in $\mathbb{F}_q^n$ vs independent sets in $\mathbb{P}^{n-1}(\mathbb{F}_q)$.} Let $2 \leq \ell \leq k \leq \frac{q^n-1}{q-1}$ with $\ell \leq n$.  Then:
\begin{enumerate}
\item ${\rm Ind}_q^{{\rm pro}}(n,k,\ell) \leq {\rm Ind}_q(n,k,\ell)$, with equality when $k=\ell$; \label{thm:mainindaffproj3a}
\item
$
        \left\lceil\frac{1}{q-1}\operatorname{Ind}_q(n,k,\ell)\right\rceil\leq \operatorname{Ind}_q^{\operatorname{pro}}(n,k,\ell)\leq \left\lfloor\frac{1}{q-1}\operatorname{Ind}_q(n,(q-1)k,\ell)\right\rfloor.
    $ Moreover, if $k\leq \frac{q^{\ell-1}-1}{q-1}$, then ${\rm Ind}_q(n,k,\ell)=(q-1){\rm Ind}^{{\rm pro}}_q(n,k,\ell)$ if and only if $q=2$. \label{thm:mainindaffproj3b}
\end{enumerate}
\end{enumerate} 
\end{theorem}
This paper is organized as follows. In Section \ref{sec:indsetaffproj}, we investigate independent sets in $\mathbb{F}_q^n$ and $\mathbb{P}^{n-1}(\mathbb{F}_q)$; in particular we prove Theorem \ref{thm:mainindaffproj}. Section \ref{sec:orthsetsub} is devoted to orthogonality in $\mathcal{K}^n$:
\begin{itemize}
\item {\bf Weakly orthogonal sets and feebly orthogonal spaces:} Theorem \ref{thm:main1} and Theorem \ref{thm:main1ss} are proved in \S\ref{sec:mainpart1}.
\item {\bf Orthogonal sets:} Theorem \ref{thm:main2} is proved in \S\ref{sec:mainpart2}.
\item {\bf Strongly orthogonal sets:} Theorem \ref{thm:main3} is proved in \S\ref{sec:mainpart3}.
\end{itemize}
\section{Independent sets in $\mathbb{F}_q^n$ and $\mathbb{P}^{n-1}(\mathbb{F}_q)$}\label{sec:indsetaffproj}
In this section, we deal with independent sets in $\mathbb{F}_q^n$ and $\mathbb{P}^{n-1}(\mathbb{F}_q)$. In particular, we prove Theorem \ref{thm:mainindaffproj}. Let $\mathbb{F}_q$ be a finite field and let $n\geq 1$ be an integer.
\vskip 1mm
\noindent
{\bf Notations:}
\begin{itemize}
\item For $m \geq 1$, we let $[m]=\{1,\dots,m\}$.
\item For ${\bf v}=(v_1,\dots,v_n) \in \mathbb{F}_q^n$, let ${\rm supp}({\bf v})=\{i\in [n] \mid v_i \neq 0\}$.
\item Let  $ s \in [n]$. The \emph{Grassmannian} ${\rm Gr}_{s,n}(\mathbb{F}_q)$ denotes the space of all $s$-dimensional subspaces of $\mathbb{F}_q^n$; the size of ${\rm Gr}_{s,n}(\mathbb{F}_q)$ is given by the \emph{$q$-binomial} ${n \brack s}_q=\frac{(q^{n-s+1}-1)\dots(q^2-1)(q-1)}{(q^s-1)\dots(q^2-1)(q-1)}$.
\item Let $\rho_n:\mathbb{F}_q^n\setminus \{{\bf 0}\}\rightarrow \mathbb{P}^{n-1}(\mathbb{F}_q)$ be the natural projection. 
\end{itemize}
\subsection{Independent sets in $\mathbb{F}_q^n$} In this part, we establish some properties of ${\rm Ind}_q(n,k,\ell)$.
In \cite{DMMS} and \cite{DMM}, Damelin, Michalski, M\"{u}llen and Stone investigated ${\rm Ind}_q(n,k,\ell)$ when $k=\ell$. Independently,  Tassa and Villar \cite{TV09} also studied ${\rm Ind}_q(n,k,\ell)$. In finite geometry, similar concepts closely related to the $(k,\ell)$-(pro-)independence have been also investigated earlier (see the survey \cite{H} of Hirschfeld); for example, given $1 \leq k \leq n$, Tallini \cite{Tal} studied the maximum size of $S \subset \mathbb{F}_q^n$ such that any $X \subset S$ of size $k$ is $\mathbb{F}_q$-linearly independent, but not all $X' \subset S$ of size $k+1$ is $\mathbb{F}_q$-linearly independent.   
\begin{theorem}[Theorem 2 in \cite{DMMS}]
\label{thm:DMMS}
    \begin{enumerate}
        \item For every $n\geq 3$, we have ${\rm Ind}_2(n,3,3)=2^{n-1}$.\label{case:3DMMS}
        \item \label{case:n+1} For every $m\geq 0$ and $n\geq 3m+2$, we have ${\rm Ind}_2(n,n-m,n-m)=n+1$.
        \item For every $m\geq 2$, $i\in \{0,1\}$ and $n=3m+i$, we have ${\rm Ind}_2(n,n-m,n-m)=n+2$.
    \end{enumerate}
\end{theorem}
\noindent
\begin{theorem}[Theorem 2 in \cite{DMM}]
\label{thm:DMM} Let $2 \leq \ell \leq n$. Then $\operatorname{Ind}_q(n,\ell,\ell)=n+1$ if and only if $\frac{q}{q+1}(n+1)\leq \ell$.
\end{theorem}
\begin{proposition}[Theorem \ref{thm:mainindaffproj} \eqref{thm:mainindaffproj1a}]\label{lem:stabilizationInd} Let $2 \leq \ell \leq k \leq q^n-1$ with $\ell \leq n$. Then $\operatorname{Ind}_q(n,k,\ell)=q^{n}-1$ if and only if $k\geq q^{\ell-1}$. In other words, $\mathbb{F}_q^{n}\setminus \{{\bf  0}\}$ is $(k,\ell)$-independent if and only if $k\geq q^{\ell-1}$.
\end{proposition}
\begin{proof}
    First assume that $k< q^{\ell-1}$. Let us consider an $\ell-1$-dimensional subspace $V \leq \mathbb{F}_q^n$. Since $\vert V \setminus \{{\bf 0}\}\vert = q^{\ell-1}-1\geq k$ and $\dim(\langle V\setminus \{\bf 0\} \rangle) < \ell$, the set $\mathbb{F}_q^n\setminus\{\boldsymbol{0}\}$ is not $(k,\ell)$-independent.

    Now assume that $k \geq q^{\ell-1}$. Let $S\subseteq \mathbb{F}_q^n\setminus\{{\bf 0}\}$ be a subset of size $k$. Since $\vert S\cup\{{\bf 0}\}\vert =k+1> q^{\ell-1}$, we have $\dim(\langle S\rangle) \geq \ell$. Hence $\mathbb{F}_q^n\setminus\{{\bf 0}\}$ is $(k,\ell)$-independent.

    This completes the proof.
\end{proof}
\begin{corollary}
For $k \geq 3$,
$$
{\rm Ind}_2(n,k,3)=
\begin{cases}
2^{n-1} &\textrm{if}\,\, k=3,\\
2^n-1  & \textrm{if}\,\, k \geq 4.
\end{cases}
$$
\end{corollary}
\begin{proof}
This follows from Theorem \ref{thm:DMMS} \eqref{case:3DMMS} and Proposition \ref{lem:stabilizationInd}.
\end{proof}
\begin{proposition}[Theorem \ref{thm:mainindaffproj} \eqref{thm:mainindaffproj1b} and \eqref{thm:mainindaffproj1c}]
\label{lem:IndInc}
\begin{enumerate}
\item \label{lem:IndIncFixl} Let $2 \leq \ell \leq n$. Then the sequence $\{{\rm Ind}_q(n,k,\ell)\}_{k=\ell}^{q^{\ell-1}}$ is strictly increasing in $k$.
\item \label{lem:IndIncFixk} Let $2 \leq k \leq q^n-1$. Then
$$
q^n-1={\rm Ind}_q(n,k,1) =\dots={\rm Ind}_q(n,k,\lfloor \log_q(k) \rfloor+1)> \dots > {\rm Ind}_q(n,k,k)\geq n+1.
$$
\end{enumerate}
\end{proposition}

\begin{proof}
    \begin{enumerate}
        \item First, note that by Proposition \ref{lem:stabilizationInd}, ${\rm Ind}_q(n,k,q^{\ell-1}-1)< {\rm Ind}_q(n,k,q^{\ell-1})=q^n-1$. Also, $\{{\rm Ind}_q(n,k,\ell)\}_{k=\ell}^{q^{\ell-1}}$ is weakly increasing. Assume on the contrary that it is not strictly increasing, that is, $\operatorname{Ind}_q(n,k,\ell)=\operatorname{Ind}_q(n,k+1,\ell)$ for some $\ell \leq k \leq q^{\ell-1}-2$. Let us consider a $(k,\ell)$-independent set $S \subseteq \mathbb{F}_q^n\setminus\{{\bf 0}\}$ of maximum size, and let $\mathbf{v}\in \mathbb{F}_q^n\setminus (S\cup\{{\bf 0}\})$.
        \begin{claim*}
        The set $S\cup\{\mathbf{v}\}$ is $(k+1,\ell)$-independent.
        \end{claim*}
        \begin{proof}[Proof of the claim]
        Let  $X\subseteq S\cup\{\mathbf{v}\}$ be of size $k+1$. Since $\vert S \cap X\vert  \geq k$, by the $(k,\ell)$-independence of $S$, the set $S \cap X$ contains an independent set of size $\ell$. Hence $S\cup\{\mathbf{v}\}$ is $(k+1,\ell)$-independent, which completes the proof of the claim.
        \end{proof}
    
        Therefore $\operatorname{Ind}_q(n,k+1,\ell)\geq \vert S\vert +1>\operatorname{Ind}_q(n,k,\ell)$, a contradiction. Consequently $\{{\rm Ind}_q(n,k,\ell)\}_{k=\ell}^{q^{\ell-1}}$ is strictly increasing.
        \item Clearly $\{{\rm Ind}_q(n,k,\ell)\}_{\ell=1}^k$ is weakly decreasing in $\ell$. By Proposition \ref{lem:stabilizationInd}, we have 
        $$q^n-1={\rm Ind}_q(n,k,1) =\dots={\rm Ind}_q(n,k,\lfloor\log_q(k)\rfloor+1)
        $$ and 
        $${\rm Ind}_q(n,k,\lfloor\log_q(k)\rfloor+1)> {\rm Ind}_q(n,k,\lfloor\log_q(k)\rfloor+2).$$
        Assume on the contrary that ${\rm Ind}_q(n,k,\ell)={\rm Ind}_q(n,k,\ell+1)$ for some $ \lfloor \log_q(k)\rfloor +2 \leq \ell \leq k-1$. Note that ${\rm Ind}_q(n,k,\ell)<q^n-1$. Now let us consider a $(k,\ell+1)$-independent set $S$ of maximum size, and let $\mathbf{v}\in \mathbb{F}_q^n\setminus (S\cup\{{\bf 0}\})$. 
        \begin{claim*}
        The set $S\cup \{\mathbf{v}\}$ is $(k,\ell)$-independent.
        \end{claim*}
        \begin{proof}[Proof of the claim]
        Let $X\subseteq S\cup \{\mathbf{v}\}$ be of size $k$. If $X \subseteq S$, then, by the $(k,\ell+1)$-independence of $S$, the set $X$ contains an independent set of size $\ell$. We may then assume that $X \not\subset S$. Let us choose ${\bf u} \in S \setminus X$, which exists because $\vert S\vert  \geq k$. Since $\vert (X\cap S) \cup \{{\bf u}\}\vert\geq k$, by the $(k,\ell+1)$-independence of $S$, the set $(X\cap S)\cup\{{\bf u}\}$ contains an independent set of size $\ell+1$. Hence $X$ contains an independent set of size $\ell$. This implies that $S\cup \{\mathbf{v}\}$ is $(k,\ell)$-independent.
        \end{proof}
        Therefore ${\rm Ind}_q(n,k,\ell)\geq\vert S\vert +1 > {\rm Ind}_q(n,k,\ell+1)$, a contradiction. Consequently $$
q^n-1={\rm Ind}_q(n,k,1) =\dots={\rm Ind}_q(n,k,\lfloor \log_q(k) \rfloor+1)> \dots > {\rm Ind}_q(n,k,k)\geq n+1.
$$
    \end{enumerate}
\end{proof}

The following generalizes slightly Theorem \ref{thm:DMM}.
\begin{corollary}Let $2 \leq \ell \leq k \leq q^{n}-1$. Then $\operatorname{Ind}_q(n,k,\ell)=n+1$ if and only if $k=\ell$ and $\frac{q}{q+1}(n+1) \leq \ell$.
\end{corollary}
\begin{proof}
    This is a consequence of Theorem \ref{thm:DMM} and Proposition \ref{lem:IndInc} \eqref{lem:IndIncFixl}.
\end{proof}
\subsection{Independent sets in $\mathbb{P}^{n-1}(\mathbb{F}_q)$} In this part, we investigate some properties of $\operatorname{Ind}_q^{\operatorname{pro}}(n,k,\ell)$.
\begin{proposition}[Theorem \ref{thm:mainindaffproj} \eqref{thm:mainindaffproj2a}] Let $2 \leq \ell \leq k \leq \frac{q^n-1}{q-1}$ with $\ell \leq n$. Then $\operatorname{Ind}_q^{\operatorname{pro}}(n,k,\ell)=\frac{q^n-1}{q-1}$ if and only if $k> \frac{q^{\ell-1}-1}{q-1}$. In other words, $\mathbb{P}^{n-1}(\mathbb{F}_q)$ is $(k,\ell)$-pro-independent if and only if $k> \frac{q^{\ell-1}-1}{q-1}$. 
\end{proposition}
\begin{proof}
    First assume that $k\leq \frac{q^{\ell-1}-1}{q-1}$. Let us consider an $\ell-1$-dimensional subspace $V \leq \mathbb{F}_q^{n}$. Since $\vert \rho_n(V \setminus \{{\bf 0}\})\vert =\frac{q^{\ell-1}-1}{q-1}\geq k$ and $\dim(\langle X \rangle)< \ell$, the set $\mathbb{P}^{n-1}(\mathbb{F}_q)$ is not $(k,\ell)$-pro-independent. 

    Now assume that $k> \frac{q^{\ell-1}-1}{q-1}$. Let $S\subseteq \mathbb{P}^{n-1}(\mathbb{F}_q)$ be of size $k$. Since $\vert S\vert =k > \frac{q^{\ell-1}-1}{q-1}$, we have $\dim(\langle S \rangle) \geq \ell$. Hence $\mathbb{P}^{n-1}(\mathbb{F}_q)$ is $(k,\ell)$-pro-independent.
\end{proof}
\begin{proposition}[Theorem \ref{thm:mainindaffproj} \eqref{thm:mainindaffproj2b} and \eqref{thm:mainindaffproj2c}]
\label{lem:IndIncpro}
\begin{enumerate}
\item \label{lem:proIndIncFixl} Let $2 \leq \ell \leq n$. Then the sequence $\{\operatorname{Ind}_q^{\operatorname{pro}}(n,k,\ell)\}_{k=\ell}^{\frac{q^{\ell-1}-1}{q-1}}$ is strictly increasing in $k$. 
\item \label{lem:proIndIncFixk} Let $2 \leq k \leq \frac{q^n-1}{q-1}$. Then
$$
\frac{q^n-1}{q-1}={\rm Ind}^{\operatorname{pro}}_q(n,k,1)= \dots={\rm Ind}^{{\rm pro}}_{q}(n,k,\lfloor \log_q((q-1)k+1)\rfloor+1)> \dots > {\rm Ind}^{\operatorname{pro}}_q(n,k,k)\geq n+1.
$$
\end{enumerate}
\end{proposition}
\begin{proof}
    \begin{enumerate}
          
        \item The proof is analogous to the one of Proposition \ref{lem:IndInc} \eqref{lem:IndIncFixl}.
        \item The proof is analogous to the one of Proposition \ref{lem:IndInc}  \eqref{lem:IndIncFixk}. 
    \end{enumerate}
\end{proof}

\subsection{$\mathbb{F}_q^n$ versus $\mathbb{P}^{n-1}(\mathbb{F}_q)$} In this part, we discuss some relations between ${\rm Ind}_q(n,k,\ell)$ and ${\rm Ind}^{{\rm pro}}_q(n,k,\ell)$. The following is straightforward.
\begin{lemma}[Properties of $\rho_n$]
\label{lem:propPi}
    \begin{enumerate}
        \item The map $\rho_n$ is a $(q-1)$-to-$1$. As a consequence, for every $S\subseteq \mathbb{P}^{n-1}(\mathbb{F}_q)$, we have that $\vert \rho_n^{-1}(S)\vert=(q-1)\vert S\vert$.
        \item \label{propPiImSize} For any $T\subseteq \mathbb{F}_q^n \setminus \{{\bf 0}\}$, we have $\frac{\vert T\vert}{q-1}\leq\vert \rho_n (T)\vert\leq \vert T\vert$.
        \item \label{indInj}If $T\subseteq \mathbb{F}_q^n\setminus\{{\bf 0}\}$ is linearly independent, then $\vert \rho_n(T)\vert=\vert T\vert$. \label{eq:lem:propPi3}
        \item Let ${\bf v}_1,\dots,{\bf v}_r \in \mathbb{F}_q^n \setminus\{{\bf 0}\}$. If $\rho_n({\bf v}_1),\dots,\rho_n({\bf v}_r)$ are linearly independent, then so are ${\bf v}_1,\dots,{\bf v}_r$. \label{eq:lem:propPi4}
    \end{enumerate}
\end{lemma}
\begin{lemma}\label{lem:ellto2ind}
Let $2 \leq \ell \leq n$. Suppose that $S\subseteq \mathbb{F}_q^n\setminus \{{\bf 0}\}$ is an $(\ell,\ell)$-independent set. Then $S$ is $(2,2)$-independent.
\end{lemma}
\begin{proof}
Let ${\bf u},{\bf v} \in S$ be distinct vectors. Since $\vert S\vert  \geq \ell$, there exists a subset $X \subseteq S$ of size $\ell$ containing both ${\bf u}$ and ${\bf v}$. Noticing that $X$ is an independent set, ${\bf u}$ and ${\bf v}$ are also independent. Therefore $S$ is $(2,2)$-independent.
\end{proof}
\begin{lemma}
\label{lem:pi(S)proInd} Let $2 \leq \ell \leq k \leq q^n-1$ with $\ell \leq n$. Suppose that $S\subseteq \mathbb{F}_q^n\setminus \{{\bf 0}\}$ is $(k,\ell)$-independent. Then $\rho_n(S) \subseteq \mathbb{P}^{n-1}(\mathbb{F}_q)$ is $(k,\ell)$-pro-independent. Moreover, if $k=\ell$, then $\vert \rho_n(S)\vert =\vert S\vert $. 
\end{lemma}
\begin{proof}
For the first statement, let $X\subseteq \rho_n(S)$ be of size $k$. Since $\rho_n^{-1}(X) \cap S \subseteq S$ has size at least $k$, it contains an independent subset $Y$ of size $\ell$. By Lemma \ref{lem:propPi} \eqref{indInj} and \eqref{eq:lem:propPi4}, $\rho_n(Y) \subseteq X$ is also independent of size $\ell$. Consequently $\rho_n(S)$ is $(k,\ell)$-pro-independent.  

For the second statement, by Lemma \ref{lem:ellto2ind}, $S$ is $(2,2)$-independent. Hence, by Lemma \ref{lem:propPi} \eqref{eq:lem:propPi3}, $\rho_n\vert _S$ is injective; so $\vert \rho_n(S)\vert =\vert S\vert $
\end{proof}
\begin{proposition}[Theorem \ref{thm:mainindaffproj} \eqref{thm:mainindaffproj3a}]
\label{lem:pi(S)maxProInd} 
 Let $2 \leq \ell \leq k \leq \frac{q^n-1}{q-1}$ with $\ell \leq n$. Then ${\rm Ind}_q^{{\rm pro}}(n,k,\ell) \leq {\rm Ind}_q(n,k,\ell)$. Moreover $\operatorname{Ind}_q^{\operatorname{pro}}(n,\ell,\ell)=\operatorname{Ind}_q(n,\ell,\ell)$.
\end{proposition}
\begin{proof}
 For the first statement, let $S \subseteq \mathbb{P}^{n-1}(\mathbb{F}_q)$ be $(k,\ell)$-pro-independent of maximum size. Let us choose a lift $S' \subseteq \mathbb{F}_q^n \setminus \{{\bf 0}\}$ of $S$ under $\rho_n$. Since $\rho_n\vert _{S'}$ is injective, it is clear that $S'$ is $(k,\ell)$-independent. Hence ${\rm Ind}_q^{{\rm pro}}(n,k,\ell)=\vert S\vert =\vert S'\vert  \leq {\rm Ind}_q(n,k,\ell)$.

 For the second statement, let $S\subseteq \mathbb{F}_q^n\setminus \{{\bf 0}\}$ be $(\ell,\ell)$-independent of maximum size. By Lemma \ref{lem:pi(S)proInd}, $\rho_n(S)$ is $(\ell,\ell)$-pro-independent and $\vert \rho_n(S)\vert =\vert S\vert $. Hence ${\rm Ind}_q(n,\ell,\ell)=\vert S\vert =\vert \rho_n(S)\vert \leq {\rm Ind}_q^{{\rm pro}}(n,\ell,\ell)$. Consequently, from the first statement, we obtain $\operatorname{Ind}_q^{\operatorname{pro}}(n,\ell,\ell)=\operatorname{Ind}_q(n,\ell,\ell)$.
\end{proof}
\begin{question}
When is $\operatorname{Ind}_q^{\operatorname{pro}}(n,k,\ell)=\operatorname{Ind}_q(n,k,\ell)$?
\end{question}
\begin{corollary}
    Let $2 \leq \ell \leq n$. Then $\operatorname{Ind}_q^{\operatorname{pro}}(n,\ell,\ell)=n+1$ if and only if $\frac{q}{q+1}(n+1)\leq \ell$.
\end{corollary}
\begin{proof}
    This is a direct consequence of Theorem \ref{thm:DMM} and Proposition \ref{lem:pi(S)maxProInd}.
\end{proof}

\begin{proposition}[Theorem \ref{thm:mainindaffproj} \eqref{thm:mainindaffproj3b}]\label{thm:336} Let $2 \leq \ell \leq k \leq \frac{q^n-1}{q-1}$ with $\ell \leq n$. Then \begin{equation}
    \label{eqn:IndProIndIneq}
        \left\lceil\frac{1}{q-1}\operatorname{Ind}_q(n,k,\ell)\right\rceil\leq \operatorname{Ind}_q^{\operatorname{pro}}(n,k,\ell)\leq \left\lfloor\frac{1}{q-1}\operatorname{Ind}_q(n,(q-1)k,\ell)\right\rfloor.
    \end{equation}
\end{proposition}
\begin{proof}    
    {\bf -Upper bound}: Let $S\subseteq \mathbb{P}^{n-1}(\mathbb{F}_q)$ be a $(k,\ell)$-pro-independent set of maximum size. Let $X\subseteq \rho_n^{-1}(S)$ be of size $(q-1)k$. Since $ \vert \rho_n(X)\vert  \geq k$ there exists an independent subset $Y \subseteq \rho_n(X)$ of size $\ell$. Let us consider a lift $T \subseteq X$ of $Y$ under $\rho_n$. By Lemma \ref{lem:propPi} \eqref{eq:lem:propPi4}, $T$ is also linearly independent. Hence $\rho_n^{-1}(S)$ is $((q-1)k,\ell)$-independent, so
    \begin{equation*}
        \operatorname{Ind}_q^{\operatorname{pro}}(n,k,\ell)=\vert S\vert=\frac{\vert \rho_n^{-1}(S)\vert}{q-1}\leq \frac{\operatorname{Ind}_q(n,(q-1)k,\ell)}{q-1}. 
    \end{equation*}
    Therefore
    $\operatorname{Ind}_q^{\operatorname{pro}}(n,k,\ell)\leq \left\lfloor\frac{\operatorname{Ind}_q(n,k(q-1),\ell)}{q-1}\right\rfloor$. 

    {\bf -Lower bound:} Let $S\subseteq \mathbb{F}_q^n\setminus \{{\bf 0}\}$ be a $(k,\ell)$-independent set of maximum size. By Lemma \ref{lem:pi(S)proInd}, $\rho_n(S)$ is $(k,\ell)$-pro-independent. Since $\vert \rho_n(S)\vert\geq \frac{\vert S\vert}{q-1}$, we have $\operatorname{Ind}_q^{\operatorname{pro}}(n,k,\ell)\geq \vert \rho_n(S)\vert \geq \frac{\operatorname{Ind}_q(n,k,\ell)}{q-1}$. Therefore $\operatorname{Ind}_q^{\operatorname{pro}}(n,k,\ell)\geq \left\lceil\frac{\operatorname{Ind}_q(n,k,\ell)}{q-1}\right\rceil$.
\end{proof}
\begin{proposition}[Theorem \ref{thm:mainindaffproj} \eqref{thm:mainindaffproj3b}] Let $2 \leq \ell \leq k \leq \frac{q^n-1}{q-1}$ with $\ell \leq n$. Assume that $k\leq \frac{q^{\ell-1}-1}{q-1}$. Then ${\rm Ind}_q(n,k,\ell)=(q-1){\rm Ind}^{{\rm pro}}_q(n,k,\ell)$ if and only if $q=2$. 
\end{proposition}
\begin{proof}

    Assume that $q=2$. Then, by Remark \ref{rmk:Ind=Indpro,q=2}, we have ${\rm Ind}_q(n,k,\ell)=(q-1){\rm Ind}^{{\rm pro}}_q(n,k,\ell)$.
       
    Conversely assume that ${\rm Ind}_q(n,k,\ell)=(q-1){\rm Ind}^{{\rm pro}}_q(n,k,\ell)$. Let us consider a $(k,\ell)$-independent set $S\subseteq \mathbb{F}_q^n\setminus \{{\bf 0}\}$ of maximum size. 
    \begin{claim*}
    The set $\rho_n(S)$ is $(k,\ell)$-pro-independent of maximum size.
    \end{claim*}
    \begin{proof}[Proof of the claim]
    By Lemma \ref{lem:pi(S)proInd}, $\rho_n(S)$ is $(k,\ell)$-pro-independent. Since ${\rm Ind}^{{\rm pro}}(n,k,\ell) \geq \vert \rho_n(S)\vert  \geq \frac{\vert S\vert }{q-1}=\frac{{\rm Ind}_q(n,k,\ell)}{q-1}$, we get $\vert \rho_n(S)\vert ={\rm Ind}_q^{{\rm pro}}(n,k,\ell)$. Hence $\rho_n(S)$ is $(k,\ell)$-pro-independent of maximum size.  
    \end{proof}
    Set $r= \left\lceil \frac{k}{q-1} \right\rceil \geq 1$.
    \begin{claim*}
    The set $\rho_n(S)$ is $(r,\ell)$-pro-independent.
    \end{claim*}
    \begin{proof}[Proof of the claim] First, note that $(q-1)\vert \rho_n(S)\vert=\vert S\vert$. Hence $\rho_n\vert_S: S \rightarrow \rho_n(S)$ is $(q-1)$ to $1$. Now let $V \subseteq \rho_n(S)$ be of size $r$. As $X=\rho_n^{-1}(V)\cap S$ has size at least $k$, by the $(k,\ell)$-independence of $S$, we have $\dim(\langle X \rangle) \geq \ell$. Therefore $\dim(\langle V \rangle) \geq \ell$. This implies that $\rho_n(S)$ is $(r,\ell)$-pro-independent. 
    \end{proof}
    Consequently, ${\rm Ind}_q^{{\rm pro}}(n,k,\ell)=|\rho_n(S)|\leq {\rm Ind}^{{\rm pro}}_q(n,r,\ell)$. By Proposition \ref{lem:IndIncpro} \eqref{lem:proIndIncFixl}, we deduce that $r=k$, so $q=2$. 
\end{proof}
\begin{question}
    When are the inequalities in Proposition \ref{thm:336} tight?
\end{question}
\section{Orthogonality in $\mathcal{K}^n$}\label{sec:orthsetsub}
In this section we prove Theorem \ref{thm:main1}, Theorem \ref{thm:main1ss}, Theorem \ref{thm:main2} and Theorem \ref{thm:main3}. For that, let $\mathcal{K}$ be a discrete valued field with valuation $\nu$, maximal ideal $\mathfrak{m}$, valuation ring $\mathcal{O}$ and {\bf finite} residue field $\kappa=\mathbb{F}_q$. Let us fix a uniformizer $\pi$ (that is $\nu(\pi)=1$) of $\mathcal{K}$ and an integer $n \geq 1$. 
\vskip 1mm
\noindent
{\bf Notations:}
\begin{itemize}
\item Denote by $\{{\bf e}_1,\dots,{\bf e}_n\}$ the standard basis for $\mathcal{K}^n$.
\item The \emph{Grassmannians} ${\rm Gr}_{s,n}(\mathcal{K})$ and ${\rm Gr}_{n,s}(\mathbb{F}_q)$ ($1 \leq s \leq n$) denote the set of all $s$-dimensional subspaces of $\mathcal{K}^n$ and $\mathbb{F}_q^n$, respectively. In particular, ${\rm Gr}_{1,n}(\mathbb{F}_q)=\mathbb{P}^{n-1}(\mathbb{F}_q)$. 
\item For $\mathbf{v}\in \mathcal{K}^n \setminus \{{\bf 0}\}$, let $\nu({\bf v})=\min_{1 \leq i \leq n} (\nu(v_i))$.
\item Denote by $\mathbb{B}_n$ the \emph{unit sphere} in $\mathcal{K}^{n}$, that is, $\mathbb{B}_n=\{{\bf v} \in \mathcal{K}^n \mid \Vert{\bf v}\Vert=1\}$ .
\item For $m \in \mathbb{Z}$ and $\lambda \in \pi^m\mathcal{O}$, write $\lambda=O(\pi^m)$ and denote by ${\rm res}_m(\lambda)$ the residue of $\lambda$ in $\pi^m\mathcal{O}/\pi^{m+1}\mathcal{O} \simeq \mathbb{F}_q$. Let $\gamma={\rm res}_0: \mathcal{O} \rightarrow \mathbb{F}_q$ (the reduction modulo $\mathfrak{m}$) and let us fix a section $\delta:\mathbb{F}_q \rightarrow \mathcal{O}$ of $\gamma$. 
\item For ${\bf v}=(v_1,\dots,v_n) \in \mathcal{O}^n$, let $\gamma_n({\bf v})=(\gamma(v_1),\dots,\gamma(v_n)) \in \mathbb{F}_q^n$. Denote by $\delta_n:\mathbb{F}_q^n \rightarrow \mathcal{O}^n$ the section of $\gamma_n$ induced by $\delta$.
\item 
 For ${\bf 0} \neq V \subseteq \mathcal{K}^n$, let $\mu_n(V)=\gamma_n(V \cap \mathbb{B}_n)\cup \{{\bf 0}\}$.
\item Denote by $\rho_n$ the natural projection $\mathbb{F}_q^n\setminus\{{\bf 0}\} \rightarrow \mathbb{P}^{n-1}(\mathbb{F}_q)$.
\item For $\ell \geq 1$, we denote by ${\rm Sym}(\ell)$ the symmetric group on $\{1,\dots,\ell\}$.
\end{itemize}
\begin{example}\label{examp:forfunct}
    Let $\mathcal{K}=\mathbb{F}_q((x^{-1}))$ be as in Example \ref{ex:dvr}. Take $\pi=x^{-1}$. Recall that $\mathcal{O}=\mathbb{F}_q[[x^{-1}]]$ and note that $\mathbb{B}_n=\mathcal{O}^n$. We have ${\rm res}_m\left(\sum_{j=m}^{\infty}a_jx^{-j}\right)=a_m$ and 
    $$\gamma_n\left(\sum_{j=0}^{\infty}a_{1j}x^{-j},\dots,\sum_{j=0}^{\infty}a_{nj}x^{-j}\right)=\left(a_{10},\dots,a_{n0}\right).$$ A section of $\gamma$ is given by the inclusion $\delta: a \in \mathbb{F}_q \mapsto a \in \mathcal{O}$.
\end{example}
\subsection{Weak and feeble orthogonality}\label{sec:mainpart1}
In this part, we provide proofs of Theorem \ref{thm:main1} (see \S\ref{ssecproofofmainthm}) and Theorem \ref{thm:main1ss} (see \S\ref{ssecproofofmainthmss}). We begin with some preliminary results. 
\begin{lemma}
\label{lem:diffM_vOrth} Let ${\bf u}, {\bf v} \in \mathbb{B}_n$. Then $\mathbf{u}$ and $\mathbf{v}$ are orthogonal if and only if $\rho_n(\gamma_n(\mathbf{u})) \neq \rho_n(\gamma_n(\mathbf{v}))$.
\end{lemma}
\begin{proof} 
Assume that $\rho_n(\gamma_n(\mathbf{u}))=\rho_n(\gamma_n(\mathbf{v}))$. Then there exists $(a,b)\in \mathbb{F}_q^2 \setminus\{{\bf 0}\}$ such that $a \gamma_n({\bf u})+b \gamma_n({\bf v})=0$. Hence, for every $i\in [n]$, we have
$$  \vert \delta(a)u_i+\delta(b)v_i\vert=\vert O(\pi)\vert<\max\{\Vert \delta(a)\mathbf{u}\Vert,\Vert \delta(b)\mathbf{v}\Vert\}.$$ This implies that $\Vert \delta(a){\bf u}+ \delta(b) {\bf v}\Vert< \max\{\Vert \delta(a)\mathbf{u}\Vert,\Vert \delta(b)\mathbf{v}\Vert\}$. Therefore, $\mathbf{u}$ and $\mathbf{v}$ are not orthogonal.

Conversely assume that $\rho_n(\gamma_n(\mathbf{u}))\neq\rho_n(\gamma_n(\mathbf{v}))$. Let $(\lambda,\beta)\in \mathcal{K}^2 \setminus\{{\bf 0}\}$.
 \begin{itemize}
 \item {\bf Case 1: $\vert \lambda\vert <\vert \beta\vert $ or $\vert \beta\vert <\vert \lambda\vert $.} By the ultrametric inequality, 
    $$\Vert \lambda\mathbf{u}+\beta\mathbf{v}\Vert=\max\{\Vert \lambda\mathbf{u}\Vert,\Vert \beta\mathbf{v}\Vert\}.$$
    \item {\bf Case 2: $\vert \lambda\vert =\vert \beta\vert $.} Since $\rho_n(\gamma_n({\bf u})) \neq \rho_n(\gamma_n({\bf v}))$, we have ${\rm res}_{\nu(\lambda)}(\lambda)\gamma_n({\bf v})+{\rm res}_{\nu(\beta)}(\beta)\gamma_n(\mathbf{u})\neq {\bf 0}$. Fix $i_0\in {\rm supp}(\gamma_n({\bf u}))\cup {\rm supp}(\gamma_n({\bf v}))$ such that $${\rm res}_{\nu(\lambda)}(\lambda)\gamma(v_{i_0})+{\rm res}_{\nu(\beta)}(\beta)\gamma(u_{i_0}) \neq 0.$$ Since 
    \begin{equation*}
        {\rm res}_{\nu(\lambda)}\left(\lambda v_{i_0}+\beta u_{i_0}\right)={\rm res}_{\nu(\lambda)}(\lambda)\gamma(v_{i_0})+{\rm res}_{\nu(\beta)}(\beta)\gamma(u_{i_0}) \neq 0,
    \end{equation*}
     we have 
    $$\vert \lambda v_{i_0}+\beta u_{i_0}\vert=q^{\nu(\lambda)}=\max\{\Vert \lambda \mathbf{u}\Vert,\Vert \beta\mathbf{v}\Vert\},$$
    and so 
    $$\Vert \lambda \mathbf{u}+\beta\mathbf{v}\Vert=\vert \lambda u_{i_0}+\beta v_{i_0}\vert=\max\{\Vert \lambda \mathbf{u}\Vert,\Vert \beta\mathbf{v}\Vert\}.$$ 
\end{itemize}
Therefore $\Vert \lambda\mathbf{u}+\beta\mathbf{v}\Vert=\max\{\Vert \lambda \mathbf{u}\Vert,\Vert\beta\mathbf{v}\Vert\}$. Consequently $\mathbf{u}$ and $\mathbf{v}$ are orthogonal.
\end{proof}
\begin{lemma}\label{lem:lftdime} Suppose that $W=\mathbb{F}_q \cdot {\bf w}_1 \oplus \dots \oplus \mathbb{F}_q \cdot {\bf w}_s \leq \mathbb{F}_q^n$ ($s \geq 1$). Let $$V=\mathcal{K}\cdot \delta_n({\bf w}_1)+\dots +\mathcal{K} \cdot \delta_n({\bf w}_s).$$ 
Then $\dim(V)=s$ and $\mu_n(V)=W$.
\end{lemma}
\begin{proof}
Assume on the contrary that $\delta_n({\bf w}_1),\dots,\delta_n({\bf w}_s)$ are linearly dependent over $\mathcal{K}$. Then $\lambda_1 \delta_n({\bf w}_1)+\dots+\lambda_s \delta_n({\bf w}_s)=0$ for some $(\lambda_1,\dots,\lambda_s)\in \mathcal{K}^s \setminus\{{\bf 0}\}$. Fix $1 \leq i_0 \leq s$ such that $\vert \lambda_{i_0}\vert =\max_{1 \leq i \leq s}\vert \lambda_i\vert $. Hence $\sum_{i \neq i_0}\gamma(\lambda_i \lambda_{i_0}^{-1}){\bf w}_{i}+{\bf w}_{i_0}=0$, a contradiction. Therefore $\delta_n({\bf w}_1),\dots,\delta_n({\bf w}_s)$ are linearly independent over $\mathcal{K}$. Consequently $\dim(V)=s$.

Now we are going to prove that $\mu_n(V)=W$.
\begin{itemize}
\item {\bf The inclusion $W \subseteq \mu_n(V)$.} Let ${\bf w} \in W$. Write ${\bf w}=\lambda_1 {\bf w}_1+\dots+\lambda_s{\bf w}_s$ where $\lambda_1,\dots,\lambda_s \in \mathbb{F}_q$. Then ${\bf w}=\gamma_n\left(\delta(\lambda_1)\delta_n({\bf w}_1)+\dots+\delta(\lambda_s)\delta_n({\bf w}_s)\right) \in \mu_n(V)$. Hence $W \subseteq \mu_n(V)$. 
\item {\bf The inclusion $\mu_n(V) \subseteq W$.} Let ${\bf a} \in \mu_n(V) \setminus \{{\bf 0}\}$. Write ${\bf a}=\gamma_n({\bf v})$ where ${\bf v} \in V \cap \mathbb{B}_n$. Moreover write ${\bf v}=\lambda_{1}\delta_n({\bf w}_1)+\dots+\lambda_{s}\delta_n({\bf w}_s)$ where  $\lambda_1,\dots,\lambda_s \in \mathcal{K}$. Fix $1 \leq i_0 \leq s$ such that $\vert \lambda_{i_0}\vert =\max_{1\leq i \leq s}\vert \lambda_i\vert $. Since $\gamma_n(\lambda_{i_0}^{-1}{\bf v})={\bf w}_{i_0}+\sum_{i \neq i_0}\gamma(\lambda_{i_0}^{-1}\lambda_i){\bf w}_i \neq {\bf 0}$, we get $\Vert\lambda_{i_0}^{-1}{\bf v}\Vert=1$ and so $\vert \lambda_{i_0}\vert =1$. Hence $\lambda_{1},\dots,\lambda_s \in \mathcal{O}$ and ${\bf a}=\gamma_n({\bf v})=\sum_{i=1}^s \gamma(\lambda_i){\bf w}_i \in W$. Therefore $\mu_n(V) \subseteq W$.
\end{itemize}
Consequently $\mu_n(V)=W$.
\end{proof}
\begin{example}
With the same notation as in Example \ref{examp:forfunct}, if ${\bf w}_1,\dots,{\bf w}_s \in \mathbb{F}_q^n$ are $\mathbb{F}_q$-linearly independent, then 
$$\mu_n(\mathcal{K}{\bf w}_1+\dots+\mathcal{K}{\bf w}_s)=\mathbb{F}_q{\bf w}_1+\dots+\mathbb{F}_q{\bf w}_s.
$$
\end{example}
\begin{lemma}\label{lem:projsubs}
Let $V \leq \mathcal{K}^n$ be a subspace of dimension $s \geq 1$. Then $\mu_n(V)$ is a subspace of $\mathbb{F}_q^n$ of dimension $s$.
\end{lemma}
\begin{proof}
\begin{claim*}
The set $\mu_n(V)$ is a subspace of $\mathbb{F}_q^n$
\end{claim*}
\begin{proof}[Proof of the claim]
Let $\lambda \in \mathbb{F}_q^*$ and let ${\bf a}, {\bf b} \in \mu_n(V)\setminus\{{\bf 0}\}$. Write ${\bf a}=\gamma_n({\bf u})$ and ${\bf b}=\gamma_n({\bf v})$ where ${\bf u},{\bf v} \in V \cap \mathbb{B}_n$. First, noting that $\delta(\lambda){\bf u} \in V \cap \mathbb{B}_n$, we have $\lambda {\bf a}=\gamma_n(\delta(\lambda) {\bf u}) \in \mu_n(V)$. Now we are going to prove that $\lambda {\bf a}+{\bf b} \in \mu_n(V)$. Without loss of generality, we may assume $\lambda {\bf a}+{\bf b} \neq {\bf 0}$. Since $\gamma_n(\delta(\lambda){\bf u}+{\bf v})=\lambda{\bf a}+{\bf b} \neq {\bf 0}$, we have $\delta(\lambda){\bf u}+{\bf v} \in V \cap \mathbb{B}_n$, and so $\lambda{\bf a}+{\bf b} \in \mu_n(V)$. Consequently, $\mu_n(V)$ is a subspace of $\mathbb{F}_q^n$.
\end{proof}
\begin{claim*}
The dimension of $\mu_n(V)$ is $s$.
\end{claim*}
\begin{proof}[Proof of the claim]
By \cite[Lemma 3.3]{KST}, $V$ admits an orthogonal basis $\{{\bf v}_1,\dots,{\bf v}_s\}$, which we may assume to be in $\mathbb{B}_n$. By Proposition \ref{cor:c(S)linInd}, $\gamma_n({\bf v}_1),\dots,\gamma_n({\bf v}_s)$ are linearly independent. Then by Lemma \ref{lem:lftdime}, $\mu_n(V)=\mathbb{F}_q \cdot \gamma_n({\bf v}_1)\oplus\dots\oplus \mathbb{F}_q \cdot \gamma_n({\bf v}_s)$. Hence  $\dim(\mu_n(V))=s$. 
\end{proof}
\end{proof}
\begin{lemma}\label{lem:equivfeeborth}
Let $s \geq 1$. Two subspaces $U,V \in {\rm Gr}_{s,n}(\mathcal{K})$ are feebly orthogonal if and only if $\mu_n(U) \neq \mu_n(V)$.
\end{lemma}
\begin{proof}
This follows from Lemma \ref{lem:diffM_vOrth}, Lemma \ref{lem:lftdime} and Lemma \ref{lem:projsubs}.
\end{proof}
    \begin{lemma}
    \label{clm:NumMaxType} Let $s\geq 1$ and let $k \geq \ell \geq 2$.
    \begin{enumerate}
    \item Let $S\subseteq \mathbb{B}_n$. Then $S$ is $(k,\ell)$-weakly orthogonal, if and only if, for any $S'\subseteq S$ of size $k$, we have $\vert (\rho_n \circ \gamma_n)(S')\vert \geq \ell$.
    \item Let $\mathcal{F}\subseteq {\rm Gr}_{s,n}(\mathcal{K})$. Then $\mathcal{F}$ is $(k,\ell)$-feebly orthogonal, if and only if, for any $\mathcal{E}\subseteq \mathcal{F}$ of size $k$, we have $\vert \mu_n(\mathcal{E})\vert  \geq \ell$.
    \end{enumerate}
    \end{lemma}
    \begin{proof}
    \begin{enumerate}
    \item This follows from the definition of $(k,\ell)$-weak orthogonality and Lemma \ref{lem:diffM_vOrth}. 
    \item This follows from the definition of $(k,\ell)$-feeble orthogonality and Lemma \ref{lem:equivfeeborth}.
    \end{enumerate}
    \end{proof}
    
    \begin{lemma}\label{lem:ret45}
    Let $s \geq 1$ and let $k \geq \ell \geq 2$.
    \begin{enumerate}
    \item Let $S \subseteq \mathbb{B}_n$. For $ {\bf w} \in \mathbb{P}^{n-1}(\mathbb{F}_q)$, let 
    $$
    t_{{\bf w},S}=\vert \{{\bf v} \in S \mid \rho_n(\gamma_n({\bf v}))={\bf w}\}\vert .
    $$ Then the following are equivalent:
    \begin{enumerate}
    \item $S$ is $(k,\ell)$-weakly orthogonal;
    \item For any $\mathcal{I} \subseteq \mathbb{P}^{n-1}(\mathbb{F}_q)$ with $\vert \mathcal{I}\vert  \leq \ell-1$, we have $\sum_{{\bf w} \in \mathcal{I}}t_{{\bf w},S} \leq k-1$.
    \end{enumerate}\label{lem:ret451}
    \item Let $\mathcal{F} \subseteq {\rm Gr}_{n,s}(\mathcal{K})$. For $ W \in {\rm Gr}_{n,s}(\mathbb{F}_q)$, let 
    $$
    t_{W,\mathcal{F}}=\vert \{V \in \mathcal{F} \mid \mu_n(V)=W\}\vert .
    $$ Then the following are equivalent:
    \begin{enumerate}
    \item $\mathcal{F}$ is $(k,\ell)$-feebly orthogonal;
    \item For any $\mathcal{R} \subseteq {\rm Gr}_{s,n}(\mathbb{F}_q)$ with $\vert \mathcal{R}\vert  \leq \ell-1$, we have $\sum_{W \in \mathcal{R}} t_{W,S} \leq k-1$.
    \end{enumerate} \label{lem:ret452} 
    \end{enumerate}
    \end{lemma}
    \begin{proof}
\begin{enumerate}
\item By Lemma \ref{clm:NumMaxType}, $S$ is $(k,\ell)$-weakly orthogonal, if and only if,
$\vert \rho_n(\gamma_n(S'))\vert \geq \ell$, for all $S'\subseteq S$ of size $k$. This is equivalent to the following: for any $\mathcal{I} \subseteq \mathbb{P}^{n-1}(\mathbb{F}_q)$ with $\vert \mathcal{I}\vert \leq \ell-1$, we have $\vert (\rho_n \circ \gamma_n)^{-1}(\mathcal{I}) \cap S\vert  \leq k-1$. In other words, for any $\mathcal{I} \subseteq \mathbb{P}^{n-1}(\mathbb{F}_q)$ with $\vert \mathcal{I}\vert  \leq \ell-1$, we have $\sum_{{\bf w} \in \mathcal{I}}t_{{\bf w},S} \leq k-1$.
\item 
By Lemma \ref{clm:NumMaxType}, $\mathcal{F}$ is $(k,\ell)$-feebly orthogonal, if and only if,
$\vert \mu_n(\mathcal{E})\vert  \geq \ell$, for all $\mathcal{E}\subseteq \mathcal{F}$ of size $k$. This is equivalent to the following: for all $\mathcal{R} \subseteq {\rm Gr}_{s,n}(\mathbb{F}_q)$ with $\vert \mathcal{R}\vert  \leq \ell-1$, we have $\vert \mu_n^{-1}(\mathcal{R}) \cap \mathcal{F}\vert  \leq k-1$. In other words, for any $\mathcal{R} \subseteq {\rm Gr}_{s,n}(\mathbb{F}_q)$ with $\vert \mathcal{R}\vert  \leq \ell-1$, we have $\sum_{W \in \mathcal{R}}t_{W,S} \leq k-1$.
\end{enumerate}
    \end{proof}
    
    \subsubsection{Proof of Theorem \ref{thm:main1}}\label{ssecproofofmainthm}
        {\bf Lower bound of \eqref{eq:boundalphakln}.}
        Take $S_0 \subseteq \mathbb{B}_n$ satisfying $\vert (\rho_n \circ \gamma_n)^{-1}({\bf z})\cap S_0\vert  = \left\lfloor \frac{k-1}{\ell-1}\right\rfloor$ for all ${\bf z} \in \mathbb{P}^{n-1}(\mathbb{F}_q)$. Note that, for any $S' \subseteq S_0$ of size $k$, we have $(\rho_n \circ \gamma_n)(S') \geq \ell$. Thus, by Lemma \ref{clm:NumMaxType}, $S_0$ is $(k,\ell)$-weakly orthogonal.  
         Therefore
        $$
       \Delta_{n}^{(k,\ell)} \geq \vert S_0\vert  = \left\lfloor \frac{k-1}{\ell-1}\right\rfloor\frac{q^n-1}{q-1}.
        $$ 
        \vskip 1mm
        \noindent
        {\bf-Upper bound of \eqref{eq:boundalphakln}.} Let $S \subseteq \mathbb{B}_n$ be a $(k,\ell)$-weakly orthogonal set of maximum size. By Lemma \ref{lem:ret45} \eqref{lem:ret451}, we have
        $$
                \vert S\vert = \sum_{{\bf w}\in \mathbb{P}^{n-1}(\mathbb{F}_q)}t_{{\bf w},S}=\frac{\sum_{\mathcal{I} \subseteq \mathbb{P}^{n-1}(\mathbb{F}_q),\vert \mathcal{I}\vert =\ell-1}\sum_{{\bf w} \in  \mathcal{I}}t_{{\bf w},S}}{\binom{\frac{q^n-1}{q-1}-1}{\ell-2}}\leq \frac{\binom{\frac{q^n-1}{q-1}}{\ell-1}(k-1)}{\binom{\frac{q^n-1}{q-1}-1}{\ell-2}}= \frac{k-1}{\ell-1} \cdot \frac{q^n-1}{q-1},
        $$      
        and so
        $$
        \Delta_{n}^{(k,\ell)}=\vert S\vert \leq \left\lfloor\frac{k-1}{\ell-1} \cdot \frac{q^n-1}{q-1}\right\rfloor.
        $$

    \vskip 1mm
    \noindent
{\bf-Equation \eqref{eq:optimform}.}
From the discussion above, it is easy to see that
$$
\Delta^{(k,\ell)}_{n}=\max\left\lbrace\left.\sum_{i=1}^{\frac{q^n-1}{q-1}}t_i\,\,\right\vert \sum_{i \in I}t_i \leq k-1,\,\,\textrm{for all}\,\, I \subseteq \left[\frac{q^n-1}{q-1}\right]\,\, \textrm{of size}\,\, \ell-1\right\rbrace.
$$
\subsubsection{Proof of Theorem \ref{thm:main1ss}}\label{ssecproofofmainthmss} The approach is quite similar to the proof of Theorem \ref{thm:main1} in \S\ref{ssecproofofmainthm}.
\vskip 1mm
\noindent
        {\bf-Lower bound of \eqref{eq:boundalphaklnss}.}
        Take $\mathcal{F}_0 \subseteq {\rm Gr}_{s,n}(\mathcal{K})$ satisfying $\vert \mu_n^{-1}(W)\cap \mathcal{F}_0\vert = \left\lfloor \frac{k-1}{\ell-1}\right\rfloor$ for all $W \in {\rm Gr}_{s,n}(\mathbb{F}_q)$. Note that, for any $\mathcal{E} \subseteq \mathcal{F}_0$ of size $k$, we have $\vert \mu_n(\mathcal{E})\vert  \geq \ell$. By Lemma \ref{clm:NumMaxType}, $\mathcal{F}_0$ is $(k,\ell)$-feebly orthogonal. Therefore
        $$
       \Delta_{n,s}^{(k,\ell)} \geq \vert \mathcal{F}_0\vert  \geq \left\lfloor \frac{k-1}{\ell-1}\right\rfloor{n \brack s}_{q}.
        $$
        \vskip 1mm
        \noindent
        {\bf-Upper bound of \eqref{eq:boundalphaklnss}.} Let $\mathcal{F} \subseteq {\rm Gr}_{s,n}(\mathcal{K})$ be a $(k,\ell)$-feebly orthogonal set of maximum size. By Lemma \ref{lem:ret45} \eqref{lem:ret452}, we have
        $$
                \vert \mathcal{F}\vert = \sum_{W \in {\rm Gr}_{s,n}(\mathbb{F}_q)} t_{W,\mathcal{F}}=\frac{\sum_{\mathcal{R} \subseteq {\rm Gr}_{s,n}(\mathbb{F}_q),\vert \mathcal{R}\vert =\ell-1}\sum_{W \in  \mathcal{R}}t_{W,\mathcal{F}}}{\binom{\vert {\rm Gr}_{s,n}(\mathbb{F}_q)\vert -1}{\ell-2}}\leq \frac{\binom{\vert {\rm Gr}_{s,n}(\mathbb{F}_q)\vert }{\ell-1}(k-1)}{\binom{\vert {\rm Gr}_{s,n}(\mathbb{F}_q)\vert -1}{\ell-2}}=\frac{k-1}{\ell-1} \cdot \vert {\rm Gr}_{s,n}(\mathbb{F}_q)\vert ,
        $$      
        and so
        $$
        \Delta_{n,s}^{(k,\ell)}=\vert \mathcal{F}\vert \leq \left\lfloor\frac{k-1}{\ell-1} \cdot {n \brack s}_{q}\right\rfloor.
        $$
    
\vskip 1mm
\noindent
{\bf -Equation \eqref{eq:optimformss}.}
From the discussion above, it is easy to see that
$$
\Delta^{(k,\ell)}_{n,s}=\max\left\lbrace\left.\sum_{i=1}^{{n \brack s }_{q}}t_i\,\,\right\vert \sum_{i \in I}t_i \leq k-1,\,\,\textrm{for all}\,\, I \subseteq \left[{n \brack s }_{q}\right]\,\, \textrm{of size}\,\, \ell-1\right\rbrace.
$$
\subsection{Orthogonal sets }\label{sec:mainpart2}
In this section, we give a proof of Theorem \ref{thm:main2} and establish another results. Our approach is based on the notion of wedge product, and we refer the reader to \cite{RW} for more details.
\begin{lemma}[Wedge product norm] Let ${\bf v}_1,\dots, {\bf v}_\ell \in \mathcal{K}^n \setminus \{{\bf 0}\}$. Then
\begin{equation}\label{eq:normwedgee}
\Vert{\bf v}_1 \wedge \dots \wedge {\bf v}_\ell\Vert=\max_{1\leq j_1<\dots<j_{\ell}\leq n}\left\vert  \det((v_{i,j_k})_{1\leq i,k \leq \ell}) \right\vert =\max_{1\leq j_1<\dots<j_{\ell}\leq n}\left\vert \sum_{\sigma\in {\rm Sym}(\ell)}\operatorname{sgn}(\sigma)\prod_{i=1}^{\ell} v_{i,j_{\sigma(i)}}\right\vert .
\end{equation} 
In particular, if $\ell=n$, then
$$
\Vert{\bf v}_1 \wedge \dots \wedge {\bf v}_n\Vert=\left\vert \det(\mathbf{v}_1,\dots ,\mathbf{v}_n)\right\vert =\left\vert  \sum_{\sigma \in {\rm Sym}(n)} {\rm sgn}(\sigma) \prod_{i=1}^n v_{i,\sigma(i)} \right\vert .
$$ 
\end{lemma}

\begin{lemma}[Hadamard's inequality]\label{lem:hadineq}
Let ${\bf v}_1,\dots, {\bf v}_\ell \in \mathcal{K}^n \setminus \{{\bf 0}\}$. Then 
$$
\Vert{\bf v}_1 \wedge \dots \wedge {\bf v}_\ell\Vert \leq \Vert{\bf v}_1\Vert\dots\Vert{\bf v}_\ell\Vert,
$$ with equality if and only if $\{{\bf v}_1,\dots, {\bf v}_\ell\}$ is orthogonal; in this case ${\bf v}_1,\dots, {\bf v}_\ell$ are linearly independent over $\mathcal{K}$.
\end{lemma}
\begin{proof}
See \cite[Lemma 2.4]{RW}.
\end{proof}
\begin{lemma}\label{lem:stronorthequiv}
 A set $\{\mathbf{v}_1,\dots, \mathbf{v}_{\ell}\} \subseteq \mathbb{B}_n$ is orthogonal if and only if
\begin{equation}
\label{eqn:hadamard}
    \max_{1\leq j_1<\dots<j_{\ell}\leq n}\left\vert \sum_{\sigma\in {\rm Sym}(\ell)}\operatorname{sgn}(\sigma)\prod_{i=1}^{\ell} v_{i,j_{\sigma(i)}}\right\vert =1.
\end{equation}
\end{lemma}
\begin{proof}
This follows from Hadamard's inequality and \eqref{eq:normwedgee}.
\end{proof}
\noindent
Lemma \ref{lem:stronorthequiv} motivates the following.
\begin{definition}
Let $\mathbf{v}_1,\dots, \mathbf{v}_{\ell} \in \mathbb{B}_n$. For any $1\leq j_1<\dots<j_{\ell}\leq n$, we define the matrix 
 $C(\mathbf{v}_1,\dots,\mathbf{v}_{\ell};j_1,\dots,j_{\ell})\in M_{\ell}(\mathbb{F}_q)$ whose the $(i,k)$-th coordinate is $\gamma(v_{i,j_k})$, that is
 $$C(\mathbf{v}_1,\dots,\mathbf{v}_{\ell};j_1,\dots,j_{\ell})=\begin{pmatrix}
    \gamma(v_{1,j_1})&\gamma(v_{1,j_2})&\dots&\gamma(v_{1,j_\ell})\\
    \gamma(v_{2,j_1})&\gamma(v_{2,j_2})&\dots&\gamma(v_{2,j_\ell})\\
    \vdots&\dots&\ddots&\vdots\\
    \gamma(v_{\ell,j_{1}})&\gamma(v_{\ell,j_{2}})&\dots&\gamma(v_{\ell,j_{\ell}})
\end{pmatrix}.$$
\end{definition}
\begin{lemma}\label{lem:translmatricc}
Let $\mathbf{v}_1,\dots, \mathbf{v}_{\ell} \in \mathbb{B}_n$. Then $\{\mathbf{v}_1,\dots ,\mathbf{v}_{\ell}\}$ is orthogonal if and only if there exist $1 \leq j_1<\dots<j_{\ell} \leq n$ such that $C(\mathbf{v}_1,\dots,\mathbf{v}_{\ell};j_1,\dots,j_{\ell})\in \operatorname{GL}_\ell(\mathbb{F}_q)$. 
\end{lemma}
\begin{proof}
  By \eqref{eqn:hadamard}, $\{\mathbf{v}_1,\dots ,\mathbf{v}_{\ell}\}$ is orthogonal if and only if 
    \begin{equation*}
    \label{eqn:det=1}
        \max_{1 \leq j_1<\dots <j_{\ell}\leq n}\left\vert \sum_{\sigma\in {\rm Sym}(\ell)}\operatorname{sgn}(\sigma)\prod_{i=1}^{\ell}v_{i,j_{\sigma(i)}}\right\vert =1,
    \end{equation*}
    if and only if there exist $1 \leq j_1<\dots<j_{\ell}\leq n$ such that
    \begin{equation*}
    \sum_{\sigma\in {\rm Sym}(\ell)}\operatorname{sgn}(\sigma)\prod_{i=1}^{\ell}
    v_{i,j_{\sigma(i)}} \in \mathcal{O}^*.
    \end{equation*} This happens if and only if there exist $1 \leq j_1<\dots<j_{\ell}\leq n$ such that
    $$
    \sum_{\sigma\in {\rm Sym}(\ell)}\operatorname{sgn}(\sigma)\prod_{i=1}^{\ell}
    \gamma(v_{i,j_{\sigma(i)}}) \neq 0,
    $$
    if and only if there exist $1 \leq j_1<\dots<j_{\ell}\leq n$ such that $\det\left(C(\mathbf{v}_1,\dots,\mathbf{v}_{\ell};j_1,\dots,j_{\ell})\right)\neq 0$.
\end{proof}
\begin{proposition}
\label{cor:c(S)linInd}
    A set $S=\{{\bf v}_1,\dots,{\bf v}_\ell\} \subseteq \mathbb{B}_n$ with $\ell \leq n$ is orthogonal if and only if the vectors $\gamma_n({\bf v_1}),\dots,\gamma_n({\bf v_\ell}) \in \mathbb{F}_q^n$ are linearly independent. Moreover $S$ is $(k,\ell)$-orthogonal if and only if for every $X\subseteq S$ of size $k$, we have $\dim(\langle \gamma_n(X)\rangle) \geq \ell$.
\end{proposition}
\begin{proof}
This follows easily from the definition of $(k,\ell)$-orthogonality and Lemma \ref{lem:translmatricc}.
\end{proof}
\begin{proposition}[Theorem \ref{thm:main2} \eqref{thm:main21}]\label{lem:boundthetadelta} Let $k \geq \ell \geq 2$. Then 
    $$ \Theta_{n}^{(k,\ell)} \leq \Delta_{n}^{(k,\ell)} \leq  \left\lfloor\frac{k-1}{\ell-1}\cdot \frac{q^n-1}{q-1}\right\rfloor .$$ 
\end{proposition}
\begin{proof} The first inequality is straightforward. The second inequality comes from Theorem \ref{thm:main1}.
\end{proof}
\begin{remark} Note that:
\begin{itemize}
 \item $\Delta_{n}^{(k,2)}=\Theta_{n}^{(k,2)}$;
 \item and that $\Theta_{n}^{(n,n)}\geq n+1$, since $S=\{\mathbf{e}_1,\dots, \mathbf{e}_n,\sum_{i=1}^n {\bf e}_i\}\subseteq \mathcal{K}^n\setminus\{{\bf 0}\}$ is $(n,n)$-orthogonal.   
\end{itemize}
\end{remark}
\begin{question}
 When is $\Theta_{n}^{(k,\ell)} = \Delta_{n}^{(k,\ell)}$?
\end{question}
\begin{lemma}\label{lem:injorthtoaffine}
Let $2 \leq \ell \leq n$ and let $S' \subset \mathbb{B}_n$ be an $(\ell,\ell)$-orthogonal set. 
\begin{enumerate}
\item $(\rho_n \circ \gamma_n)\vert _{S'}$ is injective. In particular, $\gamma_n\vert _{S'}$ is injective;
\item $S'$ is $(m,m)$-orthogonal, for all $2 \leq m \leq \ell$.
\end{enumerate}
\end{lemma}
\begin{proof} 
\begin{enumerate}
\item Let ${\bf u}, {\bf v} \in S'$ be distinct. Since $\vert S'\vert  \geq \ell$, there exists $X \subset S'$ of size $\ell$ such that ${\bf u},{\bf v} \in X$. Since $X$ is orthogonal, by Proposition \ref{cor:c(S)linInd}, $(\rho_n \circ \gamma_n)({\bf u}) \neq (\rho_n \circ \gamma_n)({\bf v})$. Hence $(\rho_n \circ \gamma_n)\vert _{S'}$ is injective.
\item Let $2 \leq m \leq \ell$. Consider a subset $X \subset S'$ of size $m$. Since $\vert S'\vert \geq \ell$, there exists $Y \subset S'$ of size $\ell$ containing $X$. By the $(\ell,\ell)$-orthogonality of $S'$, the set $Y$ is orthogonal, so is $X$. Consequently $S'$ is $(m,m)$-orthogonal.
\end{enumerate}
\end{proof}
\begin{proposition}[Theorem \ref{thm:main2} \eqref{thm:main22}]\label{lem:boundindtheta} Let $2 \leq \ell \leq k$. Then $\operatorname{Ind}_{q}(n,k,\ell)\leq\Theta_{n}^{(k,\ell)}$. Moreover ${\rm Ind}_q(n,\ell,\ell)=\Theta_n^{(\ell,\ell)}$.
\end{proposition}
\begin{proof} For the first statement, let $S \subset \mathbb{F}_q^n \setminus \{ {\bf 0}\}$ be a $(k,\ell)$-independent set of maximum size. Then, by Proposition \ref{cor:c(S)linInd}, the lift $\delta_n(S) \subset \mathbb{B}_n$ of $S$ is $(k,\ell)$-orthogonal. Hence ${\rm Ind}(n,k,\ell)=\vert S\vert =\vert \delta_n(S)\vert  \leq \Theta_n^{(k,\ell)}$.

For the second statement, let $S' \subseteq \mathbb{B}_n$ be an $(\ell,\ell)$-orthogonal set of maximum size. On the one hand, by Lemma \ref{lem:injorthtoaffine}, we have $\vert \gamma_n(S')\vert =\vert S'\vert $. On the other hand, by Proposition \ref{cor:c(S)linInd}, $\gamma_n(S')$ is $(\ell,\ell)$-independent. Hence
$$
\Theta_n^{(\ell,\ell)} \geq \vert \gamma_n(S')\vert  = \vert S'\vert ={\rm Ind}_q(n,\ell,\ell).
$$ Therefore ${\rm Ind}_q(n,\ell,\ell)=\Theta_n^{(\ell,\ell)}$
\end{proof}
\begin{question}
When is $\operatorname{Ind}_{q}(n,k,\ell)=\Theta_{n}^{(k,\ell)}$?
\end{question}
\begin{proposition}[Theorem \ref{thm:main2} \eqref{thm:main22}]
For every $2\leq k \leq q$, we have ${\rm Ind}_{q}(n,k,2)=(k-1)\frac{q^n-1}{q-1}$.
\end{proposition}
\begin{proof} Let $ 2 \leq k \leq q$. By Proposition \ref{lem:boundthetadelta} and Proposition \ref{lem:boundindtheta}, we have $ {\rm Ind}_{q}(n,k,2)\leq (k-1)\frac{q^n-1}{q-1}$. Let us consider a lift $W \subseteq \mathbb{F}_q^n \setminus \{{\bf 0}\}$ of $\mathbb{P}^{n-1}(\mathbb{F}_q)$, and let $\lambda_1,\dots,\lambda_{k-1} \in \mathbb{F}_q^*$ be distinct elements. It is easy to see that $S=\cup_{i=1}^{k-1}\lambda_{i}W \subset \mathbb{F}_q^n \setminus \{{\bf 0}\}$ is $(k,2)$-independent of size $(k-1)\frac{q^n-1}{q-1}$. Consequently, ${\rm Ind}_{q}(n,k,2)=(k-1)\frac{q^n-1}{q-1}$.
\end{proof}
The next result follows from Theorem \ref{thm:DMMS} and Proposition \ref{lem:boundindtheta}.
\begin{lemma}\label{lem:DMMS}
Assume $q=2$.
    \begin{enumerate}
        \item For every $n\geq 3$, we have $\Theta_{n}^{(3,3)}=2^{n-1}$.\label{lem:DMMS1}
        \item For every $m\geq 0$ and $n\geq 3m+2$, we have $\Theta_{n}^{(n-m,n-m)}=n+1.$
        \item For every $m\geq 2$, $i=0,1$, and $n=3m+i$, we have $\Theta_{n}^{(n-m,n-m)}=n+2$.
    \end{enumerate}
\end{lemma}
The next result follows from Theorem \ref{thm:DMM} and Proposition \ref{lem:boundindtheta}.
\begin{lemma}
\label{eqn:alpha_ll=n+1}
Let $\ell \geq 2$. Then $\Theta_{n}^{(\ell,\ell)}=n+1$ if and only if $\frac{q}{q+1}(n+1)\leq \ell$.
\end{lemma}
\begin{proposition}[Theorem \ref{thm:main2} \eqref{thm:main22}] Let $k \geq \ell \geq 2$. Assume that $\Theta_{n}^{(k,\ell)}=\operatorname{Ind}_{q}(n,k,\ell)$. Then $k\leq q^{\ell-2}+1$.
\end{proposition}
\begin{proof}
    Let $T\subseteq \mathbb{F}_q^n\setminus\{{\bf 0}\}$ be a $(k,\ell)$-independent set of maximum size. By Proposition \ref{cor:c(S)linInd}, since $\Theta_{n}^{(k,\ell)}=\operatorname{Ind}_q(n,k,\ell)$, the lift $\delta_n(T) \subseteq \mathbb{B}_n$ of $T$ is $(k,\ell)$-orthogonal of maximum size. 
    
    Fix $\mathbf{v}\in \delta_n(T)$ and ${\bf u}\in \mathcal{K}^n$ with $0<\Vert \mathbf{u}\Vert<1$. Then $\mathbf{v}+\mathbf{u}\notin \delta_n(T)$, and so $\delta_n(T)\cup \{\mathbf{v}+\mathbf{u}\}$ is not $(k,\ell)$-orthogonal. Hence, there exists $X\subseteq \delta_n(T)$ of size $k-1$, such that $X\cup \{\mathbf{v}+\mathbf{u}\}$ does not contain an orthogonal set of size $\ell$.

    Now assume on the contrary that $\mathbf{v}$ is orthogonal to some orthogonal subset $Y=\{\mathbf{w}_1,\dots,\mathbf{w}_{\ell-1}\}\subseteq X$ of size $\ell-1$. On the one hand, by Lemma \ref{lem:hadineq}, since $Y$ and ${\bf v}+{\bf u}$ are not orthogonal, we have
    \begin{equation}\label{eq:ineqYwaa}
        \Vert (\mathbf{v}+\mathbf{u})\wedge \mathbf{w}_1\wedge \dots\wedge \mathbf{w}_{\ell-1}\Vert<\Vert \mathbf{v}+ \mathbf{u}\Vert\cdot \Vert \mathbf{w}_1\Vert\cdots \Vert \mathbf{w}_{\ell-1}\Vert.
    \end{equation}
    On the other hand, again by Lemma \ref{lem:hadineq}, 
    $$\Vert \mathbf{v}\wedge \mathbf{w}_1\wedge \dots \wedge \mathbf{w}_{\ell-1}\Vert=\Vert \mathbf{v}\Vert\Vert \mathbf{w}_1\Vert \cdots\Vert \mathbf{w}_{\ell-1}\Vert,$$ 
    and
    $$\Vert \mathbf{u}\wedge \mathbf{w}_1\wedge\dots\wedge \mathbf{w}_{\ell-1}\Vert\leq\Vert {\bf u} \Vert \Vert {\bf w}_1 \Vert \dots \Vert {\bf w}_{\ell-1} \Vert <\Vert \mathbf{v}\Vert\prod_{i=1}^{\ell-1}\Vert \mathbf{w}_1\Vert;$$ so 
    \begin{equation}\label{eq:indddd}
    \Vert (\mathbf{v}+\mathbf{u})\wedge \mathbf{w}_1\wedge\cdots\wedge \mathbf{w}_{\ell-1}\Vert=\Vert \mathbf{v}\Vert\prod_{i=1}^{\ell-1}\Vert \mathbf{w}_i\Vert=\Vert \mathbf{v}+\mathbf{u}\Vert\prod_{i=1}^{\ell-1}\Vert \mathbf{w}_i\Vert.
    \end{equation} Equations \eqref{eq:ineqYwaa} and \eqref{eq:indddd} lead to a contradiction. Hence, $\mathbf{v}$ is not orthogonal to any orthogonal set $Y\subseteq X$ of size $\ell-1$. 
    
    Therefore, by Proposition \ref{cor:c(S)linInd}, $\gamma_n(\mathbf{v})$ belongs to every $\ell-1$-dimensional subset of $\gamma_n(X)$. Consequently we have $\dim \langle \gamma_n(X)\setminus \{\gamma_n(\mathbf{v})\} \rangle \leq \ell-2$ so that 
    \begin{equation*}
        k-2\leq \vert \gamma_n(X)\setminus \{\gamma_n(\mathbf{v})\}\vert\leq q^{\ell-2}-1,
    \end{equation*}
    and so $k\leq q^{\ell-2}+1$.
\end{proof}

\begin{corollary}
Let $k \geq \ell \geq 2 $. Then ${\rm Ind}(n,k,\ell) < \Theta^{(k,\ell)}_{n}$, for all $q^{\ell-2}+2\leq k \leq q^n-1$.    
\end{corollary}

\begin{proposition}[Theorem \ref{thm:main2} \eqref{thm:main23}]
\label{thm:Theta^k,lvsTheta^l,l} Let $2 \leq \ell \leq k$. Then $\Theta_{n}^{(k,\ell)}\geq \left\lfloor\frac{k-1}{\ell-1}\right\rfloor\Theta_{n}^{(\ell,\ell)}$.
\end{proposition}
\begin{proof}
    Let $S=\{\mathbf{v}_1,\dots,\mathbf{v}_{m}\}\subseteq \mathbb{B}_n$ be an $(\ell,\ell)$-orthogonal set of maximum size. Let us choose $\omega_1,\dots, \omega_{\Big\lfloor\frac{k-1}{\ell-1}\Big\rfloor}\in \gamma^{-1}(1)$ such that the $\omega_i {\bf v}_j$'s are pairwise distinct. Consider the set 
    $$Y=\left\lbrace\omega_i\mathbf{v}_j \Bigm\vert 1 \leq i \leq \left\lfloor\frac{k-1}{\ell-1}\right\rfloor,1 \leq j \leq \vert S\vert \right\rbrace \subseteq \mathbb{B}_n,$$
    of size $\left\lfloor \frac{k-1}{\ell-1}\right\rfloor\Theta_{n}^{(\ell,\ell)}$. Note that, by Lemma \ref{lem:injorthtoaffine}, the $\gamma_n({\bf v}_j)$'s are pairwise distinct. Then $\gamma_n:Y\rightarrow \gamma_n(Y)$ is a $\Big\lfloor\frac{k-1}{\ell-1}\Big\rfloor$-$1$ function. Moreover, by Proposition \ref{cor:c(S)linInd}, $\gamma_n(Y)=\gamma_n(S)$ is $(\ell,\ell)$-independent. 
    \begin{claim*}
    The set $Y$ is $(k,\ell)$-orthogonal.
    \end{claim*}
    \begin{proof}[Proof of the claim] 
    Let $X\subseteq Y$ be a subset of size $\vert X\vert =k$. Since $\vert X\vert = k > (\ell-1)\Big\lfloor \frac{k-1}{\ell-1}\Big\rfloor$, we have $\vert \gamma_n(X)\vert  \geq \ell$. Let us consider a subset $W\subseteq \gamma_n(X)$ of size $\ell$. Since $\gamma_n(Y)$ is $(\ell,\ell)$-independent, $W\subseteq \gamma_n(S)$ is linearly independent. Letting $X'$ be a lift of $W$ in $X$, again by Proposition \ref{cor:c(S)linInd}, $X'\subseteq X$ is orthogonal of size $\ell$. Therefore $Y$ is $(k,\ell)$-orthogonal.
    \end{proof}
    Consequently
    $$\Theta_{n}^{(k,\ell)}\geq\vert Y\vert =\left\lfloor\frac{k-1}{\ell-1}\right\rfloor\Theta_{n}^{(\ell,\ell)}.
    $$
\end{proof}
\begin{proposition}[Theorem \ref{thm:main2} \eqref{thm:main23}]
\label{lem:StrongOrthBound} Let $2 \leq \ell \leq k$. Then $\Theta_{n}^{(k,\ell)}\leq (k-\ell+1)\operatorname{Ind}_{q}(n,k,\ell)$.
\end{proposition}
\begin{proof}
     Let $S\subseteq \mathbb{B}_n$ be a $(k,\ell)$-orthogonal set of maximum size.  
     \begin{claim*}
     For every ${\bf z} \in \mathbb{F}_q^n\setminus \{{\bf 0}\}$, we have $\vert \gamma_n^{-1}({\bf z}) \cap S\vert  \leq k-\ell+1$.
     \end{claim*}
     \begin{proof}[Proof of the claim]
     Assume on the contrary that $\vert \gamma_n^{-1}(\mathbf{u}) \cap S\vert \geq k-\ell+2$, for some $\mathbf{u}\in \mathbb{F}_q^n\setminus \{{\bf 0}\}$. First note that $\vert S\vert  \geq k$. Let us consider $X\subseteq S\setminus \gamma_n^{-1}(\mathbf{u})$ of size at most $\ell-2$ such that $Y=X\cup (\gamma_n^{-1}(\mathbf{u}) \cap S)$ has size at least $k$. Now take $Z \subseteq Y$ of size $\ell$. Since $\vert Z \cap \gamma_n^{-1}(\mathbf{u})\vert  \geq 2$, the set $Z$ is not orthogonal. Then $Y$ does not contain an orthogonal set of size $\ell$. This contradicts the $(k,\ell)$-orthogonality of $S$. This completes the proof of the claim.
     \end{proof}
     \noindent
     Therefore we have
    \begin{equation*}
    \label{eqn:Sc(S)Bnd}
        \vert S\vert \leq (k-\ell+1)\vert \gamma_n(S)\vert .
    \end{equation*}
    By Proposition \ref{cor:c(S)linInd}, $\gamma_n(S) \subseteq \mathbb{F}_q^n \setminus \{{\bf 0}\}$ is $(k,\ell)$-independent, so $\Theta_{n}^{(k,\ell)}=\vert S\vert \leq (k-\ell+1)\vert \gamma_n(S)\vert  \leq (k-\ell+1)\operatorname{Ind}_q(n,k,\ell)$. Consequently $\Theta_{n}^{(k,\ell)}\leq (k-\ell+1)\operatorname{Ind}_{q}(n,k,\ell)$.
\end{proof}
\begin{question}
Let $\ell \geq 2$. Does the limit $\lim_{k \rightarrow \infty}\frac{\Theta^{(k,\ell)}_{n}}{k}$ exist? Is $\left(\frac{\Theta^{k,\ell}_{n}}{k}\right)_{k \geq \ell}$ increasing in $k$?    
\end{question}
\begin{proposition}[Theorem \ref{thm:main2} \eqref{thm:main24}] Let $\ell \geq 2$. Then
$$
  \frac{q^n-1}{q^{\ell-1}-1} \leq   \limsup_{k\rightarrow \infty}\frac{\Theta_{n}^{(k,\ell)}}{k} \leq q^{n}-1.
$$
\end{proposition}
\begin{proof} {\bf -Upper bound.} By Proposition \ref{lem:StrongOrthBound} and Proposition \ref{lem:stabilizationInd},
$$\limsup_{k\rightarrow\infty}\frac{\Theta_{n}^{(k,\ell)}}{k}\leq \limsup_{k\rightarrow \infty}\frac{k-\ell+1}{k}\operatorname{Ind}_q(n,k,\ell)= \limsup_{k\rightarrow\infty}\frac{k-\ell+1}{k}(q^n-1)=q^n-1.$$

\vskip 1mm
\noindent
{\bf -Lower bound.} Let us consider a $(q^{\ell-1},\ell)$-independent subset $S \subseteq \mathbb{F}_q^n\setminus\{{\bf 0}\}$ of maximum size. Fix $r \geq 1$ and let $k_r=r(q^{\ell-1}-1)+1$. For each ${\bf z} \in S$, choose $\Omega_{{\bf z}} \subseteq \gamma_n^{-1}({\bf z})$ of size $r$. Note that $\Omega=\cup_{{\bf z} \in S}\Omega_{{\bf z}}$ has size $r\vert S\vert $. 
\begin{claim*}
The set $\Omega$ is $(k_r,\ell)$-orthogonal.
\end{claim*}
\begin{proof}[Proof of the claim]
Let $Y \subseteq \Omega$ be of size $k_r$. By the pigeonhole principle, $\vert \gamma_n(Y)\vert  \geq q^{\ell-1}$. Combining the fact that $S$ is $(q^{\ell-1},\ell)$-independent with Proposition \ref{cor:c(S)linInd}, we deduce that $Y$ contains an orthogonal subset of size $\ell$. Hence $\Omega$ is $(k_r,\ell)$-orthogonal. 
\end{proof}
\noindent
Therefore, by Proposition \ref{lem:stabilizationInd}, $\Theta_{n}^{(k_r,\ell)}\geq \vert \Omega\vert =r\vert S\vert =r  \operatorname{Ind}_q(n,q^{\ell-1},\ell)=r(q^n-1)$, so 
$$\Theta_{n}^{(k_r,\ell)}\geq \frac{k_r-1}{q^{\ell-1}-1}(q^n-1),$$
which implies 
\begin{equation*}
    \frac{\Theta_{n}^{(k_r,\ell)}}{k_r}\geq \frac{k_r-1}{k_r(q^{\ell-1}-1)}(q^n-1).
\end{equation*}

Consequently we obtain
\begin{equation*}
    \limsup_{k\rightarrow \infty}\frac{\Theta_{n}^{(k,\ell)}}{k}\geq \limsup_{r\rightarrow \infty}\frac{\Theta_n^{(k_r,\ell)}}{k_r}\geq \frac{q^n-1}{q^{\ell-1}-1}.
\end{equation*}
\end{proof}

\begin {proposition}[Theorem \ref{thm:main2} \eqref{thm:main25}] Let $\ell \geq 3$. Then $\Theta_{n}^{(\ell,\ell)}<\Delta_{n}^{(\ell,\ell)}$. 
\end{proposition}
\begin{proof} Let $S\subseteq \mathbb{B}_n$ be an $(\ell,\ell)$-orthogonal set of maximum size. First note that, by Lemma \ref{lem:injorthtoaffine}, $S$ is $(3,3)$-orthogonal. Take distinct ${\bf u}, {\bf v} \in S$. By Lemma \ref{lem:injorthtoaffine}, since $(\rho_n \circ \gamma_n)({\bf u}) \neq (\rho_n \circ \gamma_n)({\bf v})$, we have $\gamma_n({\bf u}+{\bf v}) \neq 0$, so ${\bf u}+{\bf v}\in \mathbb{B}_n$. As $\gamma_n({\bf u})$, $\gamma_n({\bf v})$ and $\gamma_n({\bf u}+{\bf v})$ are linearly dependent, by Proposition \ref{cor:c(S)linInd}, ${\bf u}$, ${\bf v}$ and ${\bf u}+{\bf v}$ are not orthogonal. Hence ${\bf u}+{\bf v} \notin S$.
\begin{claim*}
The set $S \cup \{{\bf u}+{\bf v}\}$ is $(\ell,\ell)$-weakly orthogonal.
\end{claim*}
\begin{proof}[Proof of the claim]
   Let $X\subseteq S \cup \{{\bf u}+{\bf v}\}$ be of size $\ell$. If $X\subseteq S$, then it is weakly orthogonal, by the $(\ell,\ell)$-othogonality of $S$. Now assume that $X=X'\cup \{\mathbf{u}+\mathbf{v}\}$, for some $X'\subseteq S$ of size $\ell-1$. Clearly $X'$ is weakly orthogonal. Assume on the contrary that $X$ is not weakly orthogonal, that is, there is $\mathbf{w}\in X'$ such that $\mathbf{u}+\mathbf{v}$ and $\mathbf{w}$ are not orthogonal. Then, by Proposition \ref{cor:c(S)linInd}, $\gamma_n(\mathbf{u}+\mathbf{v})$ and $\gamma_n(\mathbf{w})$ are linearly dependent, so are $\gamma_n(\mathbf{u}),\gamma_n(\mathbf{v}), \gamma_n(\mathbf{w})$. Again by Proposition \ref{cor:c(S)linInd}, $\{\mathbf{u},\mathbf{v},\mathbf{w}\}\subseteq S$ is not orthogonal, a contradiction. Hence, $X$ is weakly orthogonal. Therefore $S\cup\{{\bf u}+{\bf v}\}$ is $(\ell,\ell)$-weakly orthogonal. 
\end{proof}
   Consequently, we have $\Theta_{n}^{(\ell,\ell)}<\vert S \cup \{{\bf u},{\bf v}\}\vert  \leq \Delta_{n}^{(\ell,\ell)}$.
\end{proof}
\subsection{Strongly orthogonal sets}\label{sec:mainpart3} In this part, we prove Theorem \ref{thm:main3}. 
\begin{proposition}[Theorem \ref{thm:main3}] Let $k \geq \ell \geq 2$. Then $\Gamma_{n,q}^{(k,\ell)}={\rm Ind}^{{\rm pro}}_q(n,k,\ell)$.
\end{proposition}
\begin{proof}
By Lemma \ref{lem:diffM_vOrth}, $X \subseteq \mathbb{B}_n$ is weakly orthogonal if and only if $\vert (\rho_n \circ \gamma_n)(S)\vert =\vert S\vert $. Then a subset $S \subseteq \mathbb{B}_n$ is $(k,\ell)$-strongly orthogonal, if and only if, $\vert (\rho_n \circ \gamma_n)(S)\vert =\vert S\vert $, and by Proposition \ref{cor:c(S)linInd}, for every $X\subseteq S$ of size $k$, we have $\dim(\langle \gamma_n(X) \rangle) \geq \ell$, if and only if  $(\rho_n \circ \gamma_n)(S) \subseteq \mathbb{P}^{n-1}(\mathbb{F}_q)$ is $(k,\ell)$-pro-independent of size $\vert S\vert $. This completes the proof. 
\end{proof}
\section*{Acknowledgements}
The authors thank Bogdan Chornomaz and Roy Deutch for some discussions regarding Theorem \ref{thm:main1} and Theorem \ref{thm:main1ss}. The first author would like to heartfully thank Uri Shapira for introducing her to the function field setting during the PhD. The first author is supported by the ERC grant "Dynamics on Homogeneous Spaces" (no. 754475). The second author is grateful for the support of a Technion fellowship, of an Open University of
Israel post-doctoral fellowship, and of the Israel Science Foundation (grant no. 353/21).

The authors would also like to thank the anonymous referee for his/her comments and suggestions.
  
\end{document}